\documentclass[a4paper,11pt]{amsart}
\usepackage{mathrsfs}
\usepackage[all]{xy}
\usepackage{amsmath,amssymb,amscd,bbm,amsthm,mathrsfs}
\usepackage{graphicx}
\numberwithin{equation}{section}
\newtheorem{thm}{Theorem}[section]
\newtheorem{lem}{Lemma}[section]
\newtheorem{cor}{Corollary}[section]
\newtheorem{prop}{Proposition}[section]
\newtheorem{rem}{Remark}[section]

\makeatletter
\renewcommand{\@makefntext}[1]{\parindent 1em\noindent \hb@xt@1.8em{\hss\@makefnmark}#1}
\makeatother
\begin{document}

\title[Quaternionic M\"{o}bius left gyrogroup]{Quaternionic M\"{o}bius left gyrogroup and applications}
\author{Ruiwen Wang${}^1$}

\begin{abstract}
We introduce the quaternionic M\"{o}bius addition $\oplus$ on the quaternionic unit ball $B^{4n}$ and prove that \((B^{4n}, \oplus)\) satisfies all axioms of a gyrogroup except the loop property, i.e. it is a gyrocommutative left gyrogroup. The main difficulty is that unlike the real case, we cannot use the Clifford algebra as a tool to define M\"{o}bius addition and the gyration operator. We construct the quaternionic gyration operator in terms of matrices, instead of using Clifford algebra. As in the real case, the quaternionic M\"{o}bius left gyrogroup provide algebraic foundation for quaternionic hyperbolic space. As application, we introduce the $\sigma$-translation and construct the weighted quaternionic M\"{o}bius Laplacian operators, generalizing the Weinstein operator and Geller operator in the real and complex cases, respectively, and introduce the  \(\sigma\)-convolution, which satisfies Young-type inequalities, gyroassociativity and $\sigma$-translation invariance. Our proof is based on long but elementary calculation of quaternionic matrices, instead of relatively easy calculation of Clifford algebra in the real case. We further characterize all quaternionic M\"{o}bius left gyro-subgroups, and construct quaternionic M\"{o}bius quotient left gyrogroups. The quaternionic construction of this paper restricted to the complex ball gives us the complex M\"{o}bius left gyrogroup.
\\ \par\
\end{abstract}

 \thanks{This work is partially supported by the National Natural Science Foundation of China (No. 12471080.) and the Key Program of Fujian Provincial Natural Science Foundation of China (No.2026J002007)}
\thanks{${}^1$ School of Mathematical Sciences, Xiamen University, Fujian 361005, China, Email:rwwang@stu.xmu.edu.cn}
\maketitle
\section{Introduction}
Left gyrogroups, introduced by Ungar \cite{Foguel Ungar}, are the generalization of groups with associativity replaced by gyroassociativity.
He \cite{Ungar 2003} constructed the M\"{o}bius gyrogroup by introducing the M\"{o}bius addition on the unit disc $\mathbb{D}$, $a\oplus z:=\frac{a+z}{1+\bar{a}z}$ for $a,z\in \mathbb{D}$, and the gyration operator $\mathrm{gyr}[a, b]z=(a\oplus b)z(b\oplus a)^{-1}$.
Later, Ferreira \cite{Ferreira} constructed M\"obius gyrogroups on the real Euclidean unit ball. They provide the algebraic foundation for the Poincar\'e model of hyperbolic geometry and have found applications in harmonic analysis and relativity theory.

   A nonempty set $S$ together with a binary operation, i.e. a map from $S\times S$ to $S$, is called a {\it groupoid}.
A groupoid $(G, \oplus)$ is called a {\it left gyrogroup} \cite{Foguel Ungar} if the binary operation $\oplus$ satisfies the following axioms:\\
(1) There exists at least one element $0$ satisfying $0 \oplus a = a$ for all $a \in G$;\\
(2) For each $a \in G$, there is an element $\ominus a \in G$ such that $\ominus a \oplus a = 0$;\\
(3) For any $a, b, c \in G$, there exists a unique element $\operatorname{gyr}[a, b]c \in G$ such that $\oplus$ satisfies the {\it left gyroassociative law}
\begin{equation}\begin{aligned}\label{lem gyroassociative}
a \oplus (b \oplus c) = (a \oplus b) \oplus \operatorname{gyr}[a, b]c;
\end{aligned}\end{equation}
(4) The map $\operatorname{gyr}[a, b] : G \to G$ given by $c \mapsto \operatorname{gyr}[a, b]c$ is an automorphism of the groupoid $(G, \oplus)$ (denoted by $\operatorname{gyr}[a, b] \in \operatorname{Aut}(G, \oplus)$), i.e. for $a,b,c,d\in G$, $$\mathrm{gyr}[a,b](c\oplus d)=(\mathrm{gyr}[a,b]c)\oplus(\mathrm{gyr}[a,b]d);$$\\
(5) $\operatorname{gyr}[a, -a]={\rm id}$.

  $\operatorname{gyr}[a, b]$ is called a {\it gyration operator}. A left gyrogroup $(G, \oplus)$ is called {\it gyrocommutative} if for any $a, b \in G$, we have
\begin{equation}\begin{aligned}\label{lem commu}
a \oplus b = \operatorname{gyr}[a, b](b \oplus a).
\end{aligned}\end{equation}
For any $a, b \in G$ , if the gyration operator $\operatorname{gyr}[a, b]$ satisfies the {\it left loop property}, i.e.
\begin{equation}\begin{aligned}\label{cl loop}
\operatorname{gyr}[a, b]=\operatorname{gyr}[a\oplus b, b],
\end{aligned}\end{equation}
then $(G, \oplus)$ is called a {\it gyrogroup}.

  The M\"obius addition on the real unit ball $\mathbb{B}:=\{a\in \mathbb{R}^n; |a|<1\}$ is defined as
\begin{equation}\begin{aligned}\label{oplus}
a \oplus b = (a+b)(1+\bar{a}b)^{-1},
\end{aligned}\end{equation}
for any $a, b \in \mathbb{B}$, where $\bar{a}b$ is Clifford product, while the gyration operator is defined as
\begin{equation}\begin{aligned}\label{clifford}
\mathrm{gyr}[a, b]c=(1+a\bar{b})c(1+a\bar{b})^*,
\end{aligned}\end{equation}
where $a,b,c\in\mathbb{R}^n$, $*$ is the Clifford conjugate \cite{Ferreira}. Ferreira and Ren \cite{Ferreira Ren} further constructed the M\"{o}bius gyrogroup of the ball of radius $t$ of the paravector space $\mathbb{R}\oplus V$, where $V$ is a finite-dimensional real vector space, and characterized all gyro-subgroups, M\"obius orthogonal projectors, quotient gyrogroups and fiber bundles, etc.

  On real hyperbolic spaces, Liu and Peng \cite{Liu1} defined the $\sigma$-translation by the M\"{o}bius addition and weighted multiplicative factor, to simplify the construction of weighted Laplace-Beltrami operator. They also introduced the $\sigma$-convolution associated to the $\sigma$-translation, which was used to study the generalized Helgason-Fourier transform for radial functions. With the help of the M\"{o}bius gyrogroup, Ferreira \cite{Ferreira} further constructed the generalized Helgason-Fourier transform for non-radial functions and found the eigenfunctions for weighted Laplace operators on real hyperbolic spaces. Recently, Xia and Wang \cite{Xia} constructed the M\"{o}bius addition in octonionic space and introduced the octonionic $(\sigma,t)$-translation to prove the invariance of the weighted Laplace-Beltrami operator.

  In this paper, we will construct quaternionic M\"obius left gyrogroup on the quaternionic unit ball $B^{4n}:=\{\mathbf{q}\in \mathbb{H}^n;|\mathbf{q}|<1\}$, which is also a model of quaternionic hyperbolic space \cite{sw}. As in the real case, the left quaternionic M\"obius gyrogroup provide algebraic foundation for quaternionic hyperbolic space. Recall that a {\it quaternionic M\"{o}bius transformation} is the composition of a quaternionic unitary transformation in Sp$(n)$Sp$(1)$ and $\varphi_{\mathbf{a}}$ for some $\mathbf{a}\in B^{4n}$ \cite{Wang}, where $\varphi_{\mathbf{a}}$ is an automorphism of $B^{4n}$ given by:
\begin{equation}\begin{aligned}\label{def varphia}
\varphi_{\mathbf{a}}(\mathbf{q}):=(\mathbf{a}-\mathbf{A}_{\mathbf{a}}\mathbf{q})(1-\mathbf{a}^*\mathbf{q})^{-1}.
\end{aligned}\end{equation}
Here $\mathbf{a},\mathbf{q}\in B^{4n}$ are quaternionic $(n\times 1)$-column vectors, and
\begin{equation}\begin{aligned}\label{hermitian}
\mathbf{A}_{\mathbf{a}}=(1-S_{\mathbf{a}})|\mathbf{a}|^{-2}\mathbf{a}\mathbf{a}^*+S_{\mathbf{a}}I_n,\qquad{\rm where\;} S_{\mathbf{a}}:=(1-|\mathbf{a}|^2)^{\frac12},
\end{aligned}\end{equation}
is a quaternionic $(n\times n)$-matrix, $\mathbf{a}^*=\bar{\mathbf{a}}^t$. We define the {\it quaternionic M\"{o}bius addition} on $B^{4n}$ as
\begin{equation}\begin{aligned}\label{def a oplus b}
\mathbf{a}\oplus\mathbf{q}:=\varphi_{\mathbf{a}}\left(-\mathbf{q}\right)=(\mathbf{a}+\mathbf{A}_{\mathbf{a}}\mathbf{q})(1+\mathbf{a}^*\mathbf{q})^{-1},
\end{aligned}\end{equation}
and its inverse transformation $\ominus\mathbf{q}:=-\mathbf{q}$ and $\mathbf{a}\ominus\mathbf{q}:=\mathbf{a}\oplus(-\mathbf{q})$.

 In order to construct a gyrogroup, we need the {\it gyration operator}. Note that elements of orthogonal groups can be represented by products of elements of Clifford algebra, which are used to define the gyration operator \eqref{clifford} in the real case \cite{Ferreira}. But for elements of the complex unitary group $U(n)$ or the quaternionic unitary group Sp$(n)$, there is no such representation. So unlike the real case \eqref{oplus}, quaternionic M\"{o}bius addition \eqref{def a oplus b} is only given in terms of quaternionic matrices. However, the quaternionic gyration operator is uniquely determined by the Axiom $3$ of gyrogroups, i.e. for any $\mathbf{a}, \mathbf{b}\in B^{4n}$,
\begin{equation}\begin{aligned}\label{def gyr ominus}
\text{gyr}[\mathbf{a}, \mathbf{b}]\mathbf{c}=(\ominus(\mathbf{a}\oplus\mathbf{b}))\oplus(\mathbf{a}\oplus(\mathbf{b}\oplus\mathbf{c})).
\end{aligned}\end{equation}
By tedious calculation, we obtain the following expression of {\it quaternionic gyration operator}:
\begin{equation}\begin{aligned}\label{def gyr}
\text{gyr}[\mathbf{a}, \mathbf{b}]\mathbf{c} := \mathbf{A}_{\mathbf{a}\oplus \mathbf{b}}^{-1}(\mathbf{a}\mathbf{b}^*+\mathbf{A}_{\mathbf{a}}\mathbf{A}_{\mathbf{b}})\mathbf{c}(1+\mathbf{a}^*\mathbf{b})^{-1},
\end{aligned}\end{equation}
for $\mathbf{c}\in B^{4n}$.
It is fortunate that it is a real-linear transformation, since from the formula \eqref{def gyr ominus}, we only know it to be fractional linear.

Because of the noncommutativity of the quaternionic algebra, we are not able to prove the loop property \eqref{cl loop} for the quaternionic M\"{o}bius addition.
So we only construct a left gyrogroup and do not know whether it is a gyrogroup. But we prove a special case of loop property in Proposition \ref{aboplus}, which is sufficient for our applications (cf. Proposition \ref{lem Dmu}).
\begin{thm}\label{thm gyroqp}
Groupoid $(B^{4n}, \oplus)$ is a gyrocommutative left gyrogroup.
\end{thm}
Our proof is based on long but elementary calculation of quaternionic matrices, instead of relatively easy calculation of Clifford algebra in the real case.
  With the help of the structure of left gyrogroup, we can introduce $\sigma$-translation to construct the weighted quaternionic M\"{o}bius Laplacian operator on the quaternionic unit ball. The {\it quaternionic M\"{o}bius invariant Laplacian operator} $\triangle$ \cite{Wang} is defined by
\begin{equation}\begin{aligned}
(\triangle f)(\mathbf{a})=\triangle_0 (f\circ \varphi_{\mathbf{a}})(\mathbf{0}),
\end{aligned}\end{equation}
 for $\mathbf{a}\in B^{4n}$ and $f\in C^{2}(B^{4n})$, where $\triangle_0$ is the Euclidean Laplacian on $B^{4n}$. $\triangle$ is proved quaternionic M\"{o}bius invariant in \cite[Proposition 3.1]{Wang}, i.e.
\begin{equation}\begin{aligned}\label{u in}
\triangle(f\circ\psi)(\mathbf{a})=\triangle f(\psi(\mathbf{a})),
\end{aligned}\end{equation}
for any $\mathbf{a}\in B^{4n}$ and any quaternionic M\"{o}bius transformation $\psi$, where $f\in C^2(\Omega)$.
Define the {\it $\sigma$-translation} on $B^{4n}$ as
\begin{equation}\begin{aligned}\label{def T}
T_{\mathbf{a}}f(\mathbf{q}):=\mu_{\mathbf{a}}(\mathbf{q})f\bigg((-\mathbf{a})\oplus \mathbf{q}\bigg),
\end{aligned}\end{equation}
for $\mathbf{a},\mathbf{q}\in B^{4n}$, and the {\it $\sigma$-multiplicative factor}
\begin{equation}\begin{aligned}\label{def mu}
\mu_{\mathbf{a}}(\mathbf{q}):=\left(\frac{1-|\mathbf{a}|^2}{|1-\mathbf{a}^*\mathbf{q}|^2}\right)^{\frac{1+\sigma}{2}},
\end{aligned}\end{equation}
for $\sigma\in \mathbb{R}$. We define the {\it weighted quaternionic M\"{o}bius Laplacian operator}
\begin{equation}\begin{aligned}\label{ex form}
\triangle_{\sigma}:=& 4(1-|\mathbf{q}|^2)\left( \triangle-(1+ \sigma)( R^{0'}+R^{1'})-(\sigma-1)(\sigma+1) \right),
\end{aligned}\end{equation}
where $R^{0'}$, $R^{1'}$ are the {\it quaternionic Euler operators} (cf. \eqref{tri0}).
The operator \eqref{ex form} is $\sigma$-translation invariant by the following theorem:
\begin{thm}\label{thm invariance}
For any $f\in C^2(B^{4n})$, $\mathbf{a},\mathbf{b}\in B^{4n}$ and $\sigma\in \mathbb{R}$,
\begin{equation}\begin{aligned}
\triangle_{\sigma}(T_{\mathbf{b}} f)(\mathbf{a})= T_{\mathbf{b}}(\triangle_{\sigma} f)(\mathbf{a}).
\end{aligned}\end{equation}
\end{thm}
This is the quaternionic generalization of the Weinstein operator \cite{Liu1} and the Geller operator \cite{Liu1} in real and complex cases, respectively.
We define the {\it $\sigma$-convolution} of two measurable functions $f$ and $g$ on $B^{4n}$ by
\begin{equation}\begin{aligned}\label{def convo}
(f * g)(\mathbf{p}): = \mu_{\mathbf{p}}(\mathbf{p})\int_{B^{4n}} f(\mathbf{q}) T_{\mathbf{p}} [g_{-}](\mathbf{q}) \, dV_{\sigma}(\mathbf{q}),
\end{aligned}\end{equation}
where $g_{-}(\mathbf{q}):=g(-\mathbf{q})$, $\mathbf{p}\in B^{4n}$, and $dV_{\sigma}$ is the weighted volume element
\begin{equation}\begin{aligned}\label{def dVt}
dV_{\sigma}(\mathbf{q}):=(1-|\mathbf{q}|^2)^{\sigma-2n-1}dV(\mathbf{q}).
\end{aligned}\end{equation}
We prove Young-type inequality for the $\sigma$-convolution.
\begin{thm}\label{Young}
Let \( \sigma \leq -1 \), \( 1 \leq p, q, r \leq \infty \), \( \frac{1}{p} + \frac{1}{q} = 1 + \frac{1}{r} \) and \( s = 1 - \frac{q}{r} \). For \( f \in L^p(B^{4n}, dV_{\sigma}) \), \( g \in L^q(B^{4n}, dV_{\sigma}) \), we have \( f*g \in L^r(B^{4n}, dV_{\sigma}) \) with
\[
\|f * g\|_r \leq 2^{-1-\sigma} \|\tilde{g}\|_q^{1-s} \|g\|_q^s \|f\|_p
\]
where $\tilde{g}(\mathbf{p}) := \operatorname{ess\,sup}_{\mathbf{q} \in B^{4n}} g\circ \psi^{\sharp}_{\mathbf{q}}(\mathbf{p})$ and $\psi^{\sharp}_{\mathbf{q}}:B^{4n}\rightarrow B^{4n}$ is given by $\psi_{\mathbf{q}}(\mathbf{p})=\mathrm{gyr}[\mathbf{q},\mathbf{p}]\mathbf{p}$ for any $\mathbf{p} \in B^{4n}$.
\end{thm}

The $\sigma$-convolution also satisfies the following {\it gyroassociativity}: if \( f, g, h \in L^1(B^{4n}, dV_{\sigma}) \), then
\begin{equation}\begin{aligned}\label{as eq}
(f * (g * h))(\mathbf{a}) = \bigg(\bigg( f * g_{\mathbf{a},\mathbf{q}}^{\sharp\sharp} \bigg) *_{\mathbf{q}} h \bigg)(\mathbf{a}),
\end{aligned}\end{equation}
where $g_{\mathbf{a},\mathbf{q}}^{\sharp\sharp}(\mathbf{p}):=g\bigg(\operatorname{gyr}[\mathbf{a}, -(\mathbf{q} \oplus \mathbf{p})] \operatorname{gyr}[\mathbf{q}, \mathbf{p}]\mathbf{p}\bigg)$ and $*_{\mathbf{q}}$ is the $\sigma$-convolution with respect to the variable $\mathbf{q}$.
Compared to the real case, the proof of gyroassociativity for the $\sigma$-convolution of non-radial functions is more complicated, since the expression for the quaternionic M\"{o}bius addition is only in the form of matrices, and the quaternionic algebra is noncommutative.

  We also generalize constructions of gyro-subgroups, quotient left gyrogroups and fiber bundles \cite{Ferreira,Ferreira Ren} to the quaternionic case.
The quaternionic construction of this paper restricted to the complex ball gives us the complex M\"{o}bius left gyrogroup.

  The paper is organized as follows. In Section 2, we prove that $(B^{4n}, \oplus)$ is a gyrocommutative left gyrogroup, based on the matrix formulae \eqref{def gyr} for quaternionic gyration operator.
  In Section 3, we establish various properties for $\mu_{\mathbf{a}}$ to prove Young-type inequalities and gyroassociativity for $\sigma$-convolutions.
In Section 4, we show the weighted quaternionic M\"{o}bius Laplacian operators $\triangle_{\sigma}$ to be $\sigma$-translation invariant, i.e. weighted quaternionic M\"{o}bius invariant.
In Section 5, we characterize all M\"{o}bius left gyro-subgroups and  construct quaternionic  quotient left gyrogroups. In the appendix, we prove the technical lemmas, which are used to calculate derivatives of $\mu_{\mathbf{a}}$ to prove $\sigma$-translation invariance of $\triangle_{\sigma}$.
\section{Quaternionic M\"{o}bius addition in quaternionic unit ball}
It is easy to check that, for any $\mathbf{a}\in B^{4n}$, the quaternionic M\"{o}bius transformation $\varphi_{\mathbf{a}}$ in \eqref{def varphia} is a diffeomorphism from $\overline{B^{4n}}$ to itself and is an involution, i.e. $\varphi_{\mathbf{a}}\circ \varphi_{\mathbf{a}}=id$ (cf. \cite[Proposition 2.1]{Wang}).
The following properties of matrix $\mathbf{A}_{\mathbf{a}}$ follow from definition \eqref{hermitian} directly, and will be used frequently. If there is no confusion, we denote $\mathbf{A}_{\mathbf{a}}$ simply by $\mathbf{A}$.
\begin{lem}\label{Matrix}
Let $\mathbf{a}\in B^{4n}$. Then,\\
(1) the matrix $\mathbf{A}_{\mathbf{a}}$ has two eigenspaces span $\{\mathbf{a}\}:=\{\mathbf{a}q_0;q_0\in \mathbb{H}\}$ and $\mathbf{a}^{\bot}$ with eigenvalue 1 and $S_{\mathbf{a}}$, respectively, i.e. $\mathbf{A}_{\mathbf{a}}\mathbf{a}=\mathbf{a}$ and $\mathbf{A}_{\mathbf{a}}\mathbf{v}=S_{\mathbf{a}}\mathbf{v},$ for all $\mathbf{v}\in \mathbf{a}^{\bot}$ (cf. \cite[(2.7)]{Wang}).\\
(2) $\mathbf{A}_{\mathbf{a}}$ is a hyperhermitian, i.e. $\mathbf{A}_{\mathbf{a}}^*=\mathbf{A}_{\mathbf{a}},$ and $\mathbf{A}_{\mathbf{a}}^2=\mathbf{a}\mathbf{a}^*+S_{\mathbf{a}}^2I_n.$\\
(3) The inverse of $\mathbf{A}_{\mathbf{a}}$ satisfies $\mathbf{A}_{\mathbf{a}}^{-1}\mathbf{a}=\mathbf{a}$ and $\mathbf{A}_{\mathbf{a}}^{-1}\mathbf{v}=S_{\mathbf{a}}^{-1}\mathbf{v}$ for all $\mathbf{v}\in \mathbf{a}^{\bot}$, and
\begin{equation*}\begin{aligned}
\mathbf{A}_{\mathbf{a}}^{-1}=\frac{1}{S_{\mathbf{a}}}\left(I_n-\frac{\mathbf{a}\mathbf{a}^*}{1+S_{\mathbf{a}}}\right),\qquad(\mathbf{A}_{\mathbf{a}}^2)^{-1}=&\frac{1}{S_{\mathbf{a}}^2}(I_n-\mathbf{a}\mathbf{a}^*).
\end{aligned}\end{equation*}
\end{lem}
We have the following simple lemma about the quaternionic M\"{o}bius addition.
\begin{lem}\label{lem cancellation}
For $\mathbf{a},\mathbf{b},\mathbf{c}\in B^{4n}$, we have\\
(1) the inverse property: $\ominus(\mathbf{a}\oplus\mathbf{b})=(\ominus\mathbf{a})\oplus(\ominus\mathbf{b})$;\\
(2) the left cancellation law: $\mathbf{a} \oplus ((\ominus \mathbf{a}) \oplus \mathbf{b}) = \mathbf{b}$ and $(\ominus \mathbf{a}) \oplus (\mathbf{a} \oplus \mathbf{b}) = \mathbf{b};$
\end{lem}
\begin{proof}
By definition \eqref{hermitian}, it is direct to see $\mathbf{A}_{-\mathbf{a}}=\mathbf{A}_{\mathbf{a}}$, and so
\begin{equation*}\begin{aligned}
\ominus(\mathbf{a}\oplus\mathbf{b})=&\ominus((\mathbf{a}+\mathbf{A}_{\mathbf{a}}\mathbf{b})(1+\mathbf{a}^*\mathbf{b})^{-1})\\
=&(-\mathbf{a}+\mathbf{A}_{-\mathbf{a}}(-\mathbf{b}))(1+(-\mathbf{a}^*)(-\mathbf{b}))^{-1}=(\ominus\mathbf{a})\oplus(\ominus\mathbf{b}),\\
\mathbf{a} \oplus ((\ominus \mathbf{a}) \oplus \mathbf{b})=&\mathbf{a} \oplus ((-\mathbf{a}+\mathbf{A}_{-\mathbf{a}}\mathbf{b})(1-\mathbf{a}^*\mathbf{b})^{-1})
=\mathbf{a} \oplus (-\varphi_{\mathbf{a}}(\mathbf{b}))=\varphi_{\mathbf{a}}(\varphi_{\mathbf{a}}(\mathbf{b}))=\mathbf{b}
\end{aligned}\end{equation*}
since $\varphi_{\mathbf{a}}$ is an involution. Similarly, $(\ominus \mathbf{a}) \oplus (\mathbf{a} \oplus \mathbf{b})=\mathbf{b}$.
\end{proof}
  Recall that the {\it compact symplectic group} ${\rm Sp}(n)$ consists of all quaternionic $(n\times n)$-matrices such that $AA^*=I_n$, where $*$ is the quaternionic conjugate of matrices, and ${\rm Sp}(1)$ is the group of $\mathbb{R}$-linear transformations of $\mathbb{H}^n$ given by right multiplication by unit quaternions \cite{Co}. The isometric group preserving the quaternionic inner product is ${\rm Sp}(n){\rm Sp}(1)$.
For $\mathbf{a},\mathbf{b}\in B^{4n}$,
\begin{equation}\begin{aligned}\label{module v1}
|\mathbf{a}\oplus\mathbf{b}|=|\mathbf{b}\oplus\mathbf{a}|,
\end{aligned}\end{equation}
because we can directly check that
\begin{equation}\begin{aligned}\label{module varb}
1-|\mathbf{a}\oplus\mathbf{b}|^2=1-|\varphi_{\mathbf{a}}(-\mathbf{b})|^2=\frac{(1-|\mathbf{a}|^2)(1-|\mathbf{b}|^2)}{|1+\mathbf{a}^*\mathbf{b}|^2}
\end{aligned}\end{equation}
(cf. \cite[(5.24)]{Wang}).
In order to prove the property of the quaternionic gyration operator, we need to know properties of the very important matrix $$\mathbf{M}:=\mathbf{a}\mathbf{b}^*+\mathbf{A}\mathbf{B}$$ that appear in the expression \eqref{def gyr} of the quaternionic gyration operator $\operatorname{gyr}[\mathbf{a},\mathbf{b}]$, where $\mathbf{A}=\mathbf{A}_{\mathbf{a}}$, $\mathbf{B}=\mathbf{A}_{\mathbf{b}}$.
\begin{prop}\label{prop M}
Denote $q_{\mathbf{a},\mathbf{b}}:=1+\mathbf{a}^*\mathbf{b},\qquad\mathbf{d}:=\mathbf{a}\oplus\mathbf{b}$. Then,\\
(1) $\mathbf{M}\mathbf{M}^*=(\mathbf{d}\mathbf{d}^*+S_{\mathbf{d}}^2I_n)|q_{\mathbf{a},\mathbf{b}}|^2$;\\
(2) $(\mathbf{a}\oplus \mathbf{b})^*\mathbf{M}=q_{\mathbf{b},\mathbf{a}}^*(\mathbf{b}\oplus \mathbf{a})^*$ and $\mathbf{M}^*(\mathbf{a}\oplus \mathbf{b})=(\mathbf{b}\oplus \mathbf{a})q_{\mathbf{b},\mathbf{a}}$;\\
(3) $\mathbf{M}(\mathbf{b}\oplus \mathbf{a})=(\mathbf{a}\oplus \mathbf{b})q_{\mathbf{a},\mathbf{b}}.$
\end{prop}
\begin{proof}
(1) holds because
\begin{equation}\begin{aligned}\label{ab+AB}
\mathbf{M}\mathbf{M}^*=&(\mathbf{a}\mathbf{b}^*+\mathbf{A}\mathbf{B})(\mathbf{b}\mathbf{a}^*+\mathbf{B}\mathbf{A})
=\mathbf{a}\mathbf{a}^*|\mathbf{b}|^2+\mathbf{A}\mathbf{b}\mathbf{a}^*+\mathbf{a}\mathbf{b}^*\mathbf{A}+\mathbf{A}\mathbf{B}^2\mathbf{A}\\
=&\mathbf{a}\mathbf{a}^*|\mathbf{b}|^2+\mathbf{A}\mathbf{b}\mathbf{a}^*+\mathbf{a}\mathbf{b}^*\mathbf{A}+\mathbf{A}\mathbf{b}\mathbf{b}^*\mathbf{A}+S_{\mathbf{b}}^2(\mathbf{a}\mathbf{a}^*+S_{\mathbf{a}}^2I_n)\\
=&\mathbf{a}\mathbf{a}^*+\mathbf{A}\mathbf{b}\mathbf{a}^*+\mathbf{a}\mathbf{b}^*\mathbf{A}+\mathbf{A}\mathbf{b}\mathbf{b}^*\mathbf{A}+S_{\mathbf{a}}^2S_{\mathbf{b}}^2I_n\\
=&(\mathbf{a}+\mathbf{A}\mathbf{b})(\mathbf{a}^*+\mathbf{b}^*\mathbf{A})+S_{\mathbf{a}}^2S_{\mathbf{b}}^2I_n\\
=&(\mathbf{a}+\mathbf{A}\mathbf{b})q_{\mathbf{a},\mathbf{b}}^{-1}(q_{\mathbf{a},\mathbf{b}}^*)^{-1}(\mathbf{a}^*+\mathbf{b}^*\mathbf{A})|q_{\mathbf{a},\mathbf{b}}|^2+(1-|\mathbf{a}|^2)(1-|\mathbf{b}|^2)I_n\\
=&(\mathbf{a}\oplus\mathbf{b})(\mathbf{a}\oplus\mathbf{b})^*|q_{\mathbf{a},\mathbf{b}}|^2+(1-|\mathbf{a}\oplus\mathbf{b}|^2)|q_{\mathbf{a},\mathbf{b}}|^2I_n\\
=&(\mathbf{d}\mathbf{d}^*+S_{\mathbf{d}}^2I_n)|q_{\mathbf{a},\mathbf{b}}|^2,
\end{aligned}\end{equation}
by using Lemma \ref{Matrix} (1)-(2) for matrices $\mathbf{A}$ and $\mathbf{B}$ and \eqref{module varb}.\\
(2) The second identity is the conjugate of the first one. For the first one,
\begin{equation}\begin{aligned}\label{d*M=v}
\mathbf{d}^*\mathbf{M}=&\left(\mathbf{a}\oplus\mathbf{b}\right)^*(\mathbf{a}\mathbf{b}^*+\mathbf{A}\mathbf{B})=(1+\mathbf{b}^*\mathbf{a})^{-1}\left(\mathbf{a}^*+\mathbf{b}^*\mathbf{A}^*\right)(\mathbf{a}\mathbf{b}^*+\mathbf{A}\mathbf{B})\\
=&(1+\mathbf{b}^*\mathbf{a})^{-1}\left(|\mathbf{a}|^2\mathbf{b}^*+\mathbf{b}^*\mathbf{a}\mathbf{b}^*+\mathbf{a}^*\mathbf{B}+\mathbf{b}^*\mathbf{A}^2\mathbf{B}\right)\\
=&(1+\mathbf{b}^*\mathbf{a})^{-1}\left(|\mathbf{a}|^2\mathbf{b}^*+\mathbf{b}^*\mathbf{a}\mathbf{b}^*+\mathbf{a}^*\mathbf{B}+\mathbf{b}^*\mathbf{a}\mathbf{a}^*\mathbf{B}+S_{\mathbf{a}}^2\mathbf{b}^*\right)\\
=&(1+\mathbf{b}^*\mathbf{a})^{-1}(1+\mathbf{b}^*\mathbf{a})(\mathbf{b}^*+\mathbf{a}^*\mathbf{B})=\mathbf{b}^*+\mathbf{a}^*\mathbf{B}=q_{\mathbf{b},\mathbf{a}}^*(\mathbf{b}\oplus \mathbf{a})^*,
\end{aligned}\end{equation}
by using Lemma \ref{Matrix} (1)-(2) for $\mathbf{A}$.\\
(3) By (1) and (2), we have
\begin{equation}\begin{aligned}\label{d=Mv}
\mathbf{M}(\mathbf{b}\oplus \mathbf{a})=\mathbf{M}\mathbf{M}^*\mathbf{d}q_{\mathbf{b},\mathbf{a}}^{-1}=(\mathbf{d}\mathbf{d}^*+S_{\mathbf{d}}^2I_n)\mathbf{d}|q_{\mathbf{a},\mathbf{b}}|^2q_{\mathbf{b},\mathbf{a}}^{-1}=\mathbf{d}q_{\mathbf{a},\mathbf{b}}.
\end{aligned}\end{equation}
By definition, $q_{\mathbf{a},\mathbf{b}}^*=q_{\mathbf{b},\mathbf{a}}$. The proposition is proved.
\end{proof}
\begin{prop}\label{module identity}
Let $\mathbf{a},\mathbf{b}\in B^{4n}$.
Then, we have
\begin{equation*}\begin{aligned}
\widetilde{q}_{\mathbf{a},\mathbf{b}}:=\frac{q_{\mathbf{a},\mathbf{b}}}{|q_{\mathbf{a},\mathbf{b}}|}\in {\rm Sp}(1),\qquad \mathbf{U}_{\mathbf{a},\mathbf{b}}:=\frac{1}{|q_{\mathbf{a},\mathbf{b}}|}\mathbf{A}_{\mathbf{a}\oplus\mathbf{b}}^{-1}(\mathbf{a}\mathbf{b}^*+\mathbf{A}\mathbf{B})\in{\rm Sp}(n).
\end{aligned}\end{equation*}
Moreover, $\mathbf{U}_{\mathbf{b},\mathbf{a}}=\mathbf{U}_{\mathbf{a},\mathbf{b}}^*$ and $\widetilde{q}_{\mathbf{b},\mathbf{a}}=\widetilde{q}_{\mathbf{a},\mathbf{b}}^*$.
\end{prop}
\begin{proof}
$\widetilde{q}_{\mathbf{a},\mathbf{b}}\in {\rm Sp}(1)$ follows from $|\widetilde{q}_{\mathbf{a},\mathbf{b}}|=1$. It follows from the definition of $\mathbf{U}_{\mathbf{a},\mathbf{b}}$ that
\begin{equation}\begin{aligned}
\mathbf{U}_{\mathbf{a},\mathbf{b}}\mathbf{U}_{\mathbf{a},\mathbf{b}}^*=&|q_{\mathbf{a},\mathbf{b}}|^{-2}\mathbf{A}_{\mathbf{d}}^{-1}\mathbf{M}\mathbf{M}^*\mathbf{A}_{\mathbf{d}}^{-1}=\mathbf{A}_{\mathbf{d}}^{-1}(\mathbf{d}\mathbf{d}^*+S_{\mathbf{d}}^2I_n)\mathbf{A}_{\mathbf{d}}^{-1}\\
=&\mathbf{d}\mathbf{d}^*+(I_n-\mathbf{d}\mathbf{d}^*)=I_n,
\end{aligned}\end{equation}
by Proposition \ref{prop M} (1) and Lemma \ref{Matrix} (3).
Thus, $\mathbf{U}_{\mathbf{a},\mathbf{b}}\in{\rm Sp}(n).$ By definition, $q_{\mathbf{a},\mathbf{b}}^*=q_{\mathbf{b},\mathbf{a}}$. Hence $\widetilde{q}_{\mathbf{a},\mathbf{b}}^*=\widetilde{q}_{\mathbf{b},\mathbf{a}}$. Note that
\begin{equation}\begin{aligned}
|q_{\mathbf{b},\mathbf{a}}|\mathbf{U}_{\mathbf{b},\mathbf{a}}=&\mathbf{A}_{\mathbf{b}\oplus \mathbf{a}}^{-1}(\mathbf{b}\mathbf{a}^*+\mathbf{B}\mathbf{A})=\mathbf{A}_{\mathbf{b}\oplus \mathbf{a}}^{-1}\mathbf{M}^*\\
=&\frac{1}{S_{\mathbf{b}\oplus \mathbf{a}}}\left(I_n-\frac{(\mathbf{b}\oplus \mathbf{a})(\mathbf{b}\oplus \mathbf{a})^*}{1+S_{\mathbf{b}\oplus \mathbf{a}}}\right)\mathbf{M}^*\\
=&\mathbf{M}^*\frac{1}{S_{\mathbf{a}\oplus \mathbf{b}}}\left(I_n-\frac{\mathbf{d}\mathbf{d}^*}{1+S_{\mathbf{a}\oplus \mathbf{b}}}\right)=\mathbf{M}^*\mathbf{A}_{\mathbf{d}}^{-1}=|q_{\mathbf{a},\mathbf{b}}|\mathbf{U}_{\mathbf{a},\mathbf{b}}^*,\\
\end{aligned}\end{equation}
by Lemma \ref{Matrix} (1)-(3) and $S_{\mathbf{a}\oplus \mathbf{b}}=S_{\mathbf{b}\oplus \mathbf{a}}$ by \eqref{module v1}, and
\begin{equation*}\begin{aligned}
(\mathbf{b}\oplus \mathbf{a})(\mathbf{b}\oplus \mathbf{a})^*\mathbf{M}^*
=(\mathbf{b}\oplus \mathbf{a})q_{\mathbf{b},\mathbf{a}}\mathbf{d}^*=\mathbf{M}^*\mathbf{d}\mathbf{d}^*,
\end{aligned}\end{equation*}
where we have used Proposition \ref{prop M} (2)-(3). Hence, $\mathbf{U}_{\mathbf{b},\mathbf{a}}=\mathbf{U}_{\mathbf{a},\mathbf{b}}^*$.
\end{proof}
Now the quaternionic gyration operator can be rewritten as
\begin{equation}\begin{aligned}\label{def gyr uq}
\mathrm{gyr}[\mathbf{a}, \mathbf{b}]\mathbf{c}=\mathbf{U}_{\mathbf{a},\mathbf{b}}\mathbf{c}\widetilde{q}_{\mathbf{a},\mathbf{b}}^*.
\end{aligned}\end{equation}
So $\operatorname{gyr}[\mathbf{a}, \mathbf{b}]$ is an automorphism with respect to $\oplus$ by following lemma.
\begin{lem}\label{rotation}
For any \( \mathbf{U}\in {\rm Sp}(n),q_0\in{\rm Sp}(1)\) and \( \mathbf{a}, \mathbf{b} \in B^{4n} \), we have
\begin{equation}
\begin{aligned}
\mathbf{U}(\mathbf{a} \oplus \mathbf{b})q_0=(\mathbf{U}\mathbf{a}q_0) \oplus (\mathbf{U}\mathbf{b}q_0).
\end{aligned}
\end{equation}
\end{lem}
\begin{proof}
We denote $\widetilde{\mathbf{A}}=\mathbf{A}_{\mathbf{U}\mathbf{a}q_0}$ and $\widetilde{S}=S_{\mathbf{U}\mathbf{a}q_0}$.
Since
$$\mathbf{U}\mathbf{a}q_0(\mathbf{U}\mathbf{a}q_0)^*=\mathbf{U}\mathbf{a}\mathbf{a}^*\mathbf{U}^*,\qquad\widetilde{S}=(1-|\mathbf{U}\mathbf{a}q_0|^2)^{\frac12}=(1-|\mathbf{a}|^2)^{\frac12}=S_{\mathbf{a}},$$ we have
\begin{equation*}
\begin{aligned}
\widetilde{\mathbf{A}}=&(1-\widetilde{S})|\mathbf{U}\mathbf{a}q_0|^{-2}\mathbf{U}\mathbf{a}q_0(\mathbf{U}\mathbf{a}^*q_0)^* + \widetilde{S}I_n\\
=&\mathbf{U}\bigg(1-S_{\mathbf{a}})|\mathbf{a}|^{-2}\mathbf{a}\mathbf{a}^* + S_{\mathbf{a}}I_n\Bigg)\mathbf{U}^*=\mathbf{U}\mathbf{A}_{\mathbf{a}}\mathbf{U}^*.
\end{aligned}
\end{equation*}
Therefore,
\begin{equation}
\begin{aligned}
(\mathbf{U}\mathbf{a}q_0) \oplus (\mathbf{U}\mathbf{b}q_0)
=& \Bigl(\mathbf{U}\mathbf{a}q_0 + \widetilde{\mathbf{A}}\mathbf{U}\mathbf{b}q_0\Bigr) \cdot \bigl(1+(\mathbf{U}\mathbf{a}q_0)^*\mathbf{U}\mathbf{b}q_0\bigr)^{-1} \\
=& \mathbf{U}\Bigl(\mathbf{a}q_0 +\mathbf{A}_{\mathbf{a}}\mathbf{b}q_0\Bigr) \bigl(1+q_0^*\mathbf{a}^*\mathbf{b}q_0\bigr)^{-1}=\mathbf{U}(\mathbf{a} \oplus \mathbf{b})q_0.
\end{aligned}
\end{equation}
The lemma is proved.
\end{proof}
\begin{proof}[Proof of Theorem \ref{thm gyroqp}]
We use notations in Proposition \ref{prop M}. (1)-(2) follows from definition.\\
(3)  $\mathrm{gyr}[\mathbf{a}, \mathbf{b}]$ is a gyration operator because
\begin{equation}\begin{aligned}\label{def D}
\mathbf{a}\oplus(\mathbf{b}\oplus\mathbf{c})=&\bigg(\mathbf{a}+\mathbf{A}(\mathbf{b}+\mathbf{B}\mathbf{c})(1+\mathbf{b}^*\mathbf{c})^{-1}\bigg)\bigg(1+\mathbf{a}^*(\mathbf{b}+\mathbf{B}\mathbf{c})(1+\mathbf{b}^*\mathbf{c})^{-1}\bigg)^{-1}\\
=&\bigg(\mathbf{a}+\mathbf{A}\mathbf{b}+(\mathbf{a}\mathbf{b}^*+\mathbf{A}\mathbf{B})\mathbf{c}\bigg)\bigg(1+\mathbf{a}^*\mathbf{b}+(\mathbf{b}^*+\mathbf{a}^*\mathbf{B})\mathbf{c}\bigg)^{-1}\\
=&\bigg((\mathbf{a} + \mathbf{A}\mathbf{b})(1 + \mathbf{a}^*\mathbf{b})^{-1} + (\mathbf{a}\mathbf{b}^* + \mathbf{A}\mathbf{B})\mathbf{c}(1 + \mathbf{a}^*\mathbf{b})^{-1}\bigg)\mathfrak{F}^{-1} \\
=&\left(\mathbf{d}  + \mathbf{A}_{\mathbf{d}}\left(\mathbf{A}_{\mathbf{d}}^{-1}\mathbf{M} \mathbf{c} q_{\mathbf{a},\mathbf{b}}^{-1}\right)\right)\left(1 + \mathbf{d}^*\left( \mathbf{A}_{\mathbf{d}}^{-1}\mathbf{M} \mathbf{c} q_{\mathbf{a},\mathbf{b}}^{-1}\right)\right)^{-1}\\
=&\left(\mathbf{d}  + \mathbf{A}_{\mathbf{d}}\mathrm{gyr}[\mathbf{a}, \mathbf{b}]\mathbf{c}\right)\left(1 + \mathbf{d}^*\mathrm{gyr}[\mathbf{a}, \mathbf{b}]\mathbf{c}\right)^{-1}=(\mathbf{a} \oplus \mathbf{b}) \oplus \mathrm{gyr}[\mathbf{a}, \mathbf{b}]\mathbf{c},
\end{aligned}\end{equation}
where
\begin{equation*}\begin{aligned}
\mathfrak{F} :=& 1 + (\mathbf{b}^* + \mathbf{a}^*\mathbf{B})\mathbf{c}(1 + \mathbf{a}^*\mathbf{b})^{-1}=1 + q_{\mathbf{b},\mathbf{a}}^*(\mathbf{b}\oplus\mathbf{a})^*\mathbf{c}q_{\mathbf{a},\mathbf{b}}^{-1}\\
=& 1 + \mathbf{d}^* \mathbf{M} \mathbf{c} q_{\mathbf{a},\mathbf{b}}^{-1} = 1 + \mathbf{d}^* \mathbf{A}_{\mathbf{d}}^{-1}\mathbf{M} \mathbf{c} q_{\mathbf{a},\mathbf{b}}^{-1},
\end{aligned}\end{equation*}
by Proposition \ref{prop M} (2) and Lemma \ref{Matrix} (3) for $\mathbf{A}_{\mathbf{d}}^{-1}$.\\
  (4) $\mathrm{gyr}[\mathbf{a}, \mathbf{b}]$ is a gyroautomorphism by Lemma \ref{rotation}.\\
  (5) $(B^{4n}, \oplus)$ is gyrocommutative, since
\begin{equation*}\begin{aligned}\label{eq commu}
\mathrm{gyr}[\mathbf{a}, \mathbf{b}](\mathbf{b} \oplus \mathbf{a})=\mathbf{A}_{\mathbf{d}}^{-1}\mathbf{M}(\mathbf{b} \oplus \mathbf{a})q_{\mathbf{a},\mathbf{b}}^{-1}=\mathbf{A}_{\mathbf{d}}^{-1}\mathbf{d}=\mathbf{d}=\mathbf{a} \oplus \mathbf{b},
\end{aligned}\end{equation*}
by Lemma \ref{Matrix} (3), Proposition \ref{prop M} (3), and $q_{\mathbf{a},\mathbf{b}}=q_{\mathbf{b},\mathbf{a}}^*$ by Proposition \ref{module identity}.
\end{proof}
The quaternionic M\"{o}bius addition also satisfies following properties.
\begin{lem}\label{cancellation}
For $\mathbf{a},\mathbf{b},\mathbf{c}\in B^{4n}$, we have\\
(i) $\mathbf{a}\oplus(\mathbf{b}\oplus\mathbf{c})=\mathrm{gyr}[\mathbf{a},\mathbf{b}]((\mathbf{b}\oplus\mathbf{a})\oplus\mathbf{c});$\\
(ii) $\operatorname{gyr}[\mathbf{a}, \mathbf{b}] \operatorname{gyr}[\mathbf{b}, \mathbf{a}] = id;$\\
(iii) $(\mathbf{a}\oplus\mathbf{b})\ominus\mathbf{a}=\mathrm{gyr}[\mathbf{a},\mathbf{b}]\mathbf{b}$.
\end{lem}
\begin{proof}
(i) We have
\begin{equation}\begin{aligned}\label{2.21}
\mathbf{a}\oplus(\mathbf{b}\oplus\mathbf{c})
=&(\mathbf{a}\oplus\mathbf{b})\oplus\mathrm{gyr}[\mathbf{a},\mathbf{b}]\mathbf{c}=\mathrm{gyr}[\mathbf{a},\mathbf{b}](\mathbf{b}\oplus\mathbf{a})\oplus\mathrm{gyr}[\mathbf{a},\mathbf{b}]\mathbf{c}\\
=&\mathrm{gyr}[\mathbf{a},\mathbf{b}]((\mathbf{b}\oplus\mathbf{a})\oplus\mathbf{c}).
\end{aligned}\end{equation}
by using gyroassociativity, gyrocommutativity and gyroautomorphism, respectively.\\
(ii) For any $\mathbf{c}\in B^{4n}$, we have
$\operatorname{gyr}[\mathbf{a}, \mathbf{b}] \operatorname{gyr}[\mathbf{b}, \mathbf{a}]\mathbf{c} =\mathbf{U}_{\mathbf{a},\mathbf{b}}\mathbf{U}_{\mathbf{b},\mathbf{a}}\mathbf{c}\widetilde{q}_{\mathbf{b},\mathbf{a}}^{-1}\widetilde{q}_{\mathbf{a},\mathbf{b}}^{-1}=\mathbf{c}$ by Proposition \ref{module identity}.
(iii) By using (ii) and the identity (i) with $\mathbf{a},\mathbf{b}$ exchanged, we get
\begin{equation}\begin{aligned}
(\mathbf{a}\oplus\mathbf{b})\oplus\mathbf{c}=&\mathrm{gyr}[\mathbf{a},\mathbf{b}]\circ\mathrm{gyr}[\mathbf{b},\mathbf{a}]((\mathbf{a}\oplus\mathbf{b})\oplus\mathbf{c})\\
=&\mathrm{gyr}[\mathbf{a},\mathbf{b}](\mathbf{b}\oplus(\mathbf{a}\oplus\mathbf{c})).
\end{aligned}\end{equation}
Letting $\mathbf{c}=-\mathbf{a}$, we get (3).
\end{proof}
 The gyrocommutative gyrogroup has the following algebraic characterization. A groupoid \((L, \oplus)\) is a gyrocommutative gyrogroup if and only if it is a {\it K-loop} \cite[Theorem 3.1]{abe}. Namely, it satisfies\\
(1) For any \(a, b \in L\), the equation \(a \oplus x = b\) has a unique solution \(x \in L\);\\
(2) For any \(a, b \in L\), the equation \(y \oplus a = b\) has a unique solution \(y \in L\);\\
(3) \(a \oplus (b \oplus (a \oplus c)) = (a \oplus (b \oplus a)) \oplus c\) for any \(a, b, c \in L\);\\
(4) Any element \(a\) of \(L\) has the inverse \(\ominus a\) satisfying $\ominus(a \oplus b) = (\ominus a) \oplus (\ominus b)$  for any \(a, b \in L\).\\
We are not able to prove the loop property for the quaternionic M\"{o}bius addition, but except in the following special case, which is sufficient for the proof of propositions of the $\sigma$-translation.
\begin{prop}\label{aboplus}
For \(\mathbf{c},\mathbf{d}\in B^{4n}\), we have
\[
\operatorname{gyr}[\mathbf{c},\mathbf{d}]\mathbf{d}
=
\operatorname{gyr}[\mathbf{c}\oplus\mathbf{d},-\mathbf{c}]\mathbf{d}.
\]
\end{prop}
\begin{proof}
Applying Lemma \ref{cancellation} (iii) to $\mathbf{c},\mathbf{d}\in B^{4n}$, we get
\begin{equation}\begin{aligned}
\operatorname{gyr}[\mathbf{c},\mathbf{d}]\mathbf{d}=&(\mathbf{c}\oplus\mathbf{d})\oplus(-\mathbf{c})=\operatorname{gyr}[\mathbf{c}\oplus\mathbf{d},-\mathbf{c}]((-\mathbf{c})\oplus(\mathbf{c}\oplus\mathbf{d}))=\operatorname{gyr}[\mathbf{c}\oplus\mathbf{d},-\mathbf{c}]\mathbf{d}
\end{aligned}\end{equation}
by gyrocommutativity and left cancellation law in Lemma \ref{lem cancellation} (2). The proposition is proved.
\end{proof}
\section{$\sigma$-convolution in quaternionic M\"{o}bius left gyrogroup}
\subsection{$\sigma$-translation $T_{\mathbf{a}}$ and multiplicative factor $\mu_{\mathbf{a}}$ }
Applying \eqref{module varb} to the definition \eqref{def mu} of $\mu_{\mathbf{a}}$, we get the following useful expression
\begin{equation}\begin{aligned}\label{def mu2}
\mu_{\mathbf{a}}(\mathbf{q})=\left(\frac{1-|\varphi_{\mathbf{a}}(\mathbf{q})|^2}{1-|\mathbf{q}|^2}\right)^{\lambda}.
\end{aligned}\end{equation}
The following properties are similar to those in \cite[Lemma 1]{Ferreira} in the real case.
\begin{prop}\label{lem mu}
For $\mathbf{a},\mathbf{b},\mathbf{q}\in B^{4n}$, we have\\
(1) $\mu_{\mathbf{a}}(\mathbf{q})=\mu_{-\mathbf{a}}(-\mathbf{q}),$\\
(2) $\mu_{\mathbf{a}}(\mathbf{a})\mu_{\mathbf{a}}(\mathbf{0})=1,$\\
(3) $\mu_{\mathbf{a}}(\mathbf{q})=\mu_{\mathbf{q}}(\mathbf{a})\mu_{\mathbf{a}}(\mathbf{0})\mu_{\mathbf{q}}(\mathbf{q}),$\\
(4) $\mu_{\mathbf{a}}(\mathbf{a}\oplus\mathbf{q})=(\mu_{-\mathbf{a}}(\mathbf{q}))^{-1},$\\
(5) $\mu_{(-\mathbf{a})\oplus\mathbf{q}}(\mathbf{0})=\mu_{\mathbf{q}\oplus(-\mathbf{a})}(\mathbf{0})=\mu_{\mathbf{q}}(\mathbf{a})\mu_{\mathbf{a}}(\mathbf{0})=\mu_{\mathbf{a}}(\mathbf{q})\mu_{\mathbf{q}}(\mathbf{0}),$\\
(6) $\mu_{(-\mathbf{a})\oplus\mathbf{q}}((-\mathbf{a})\oplus\mathbf{q})=(\mu_{\mathbf{a}}(\mathbf{q}))^{-1}\mu_{\mathbf{q}}(\mathbf{q}),$\\
(7) $T_{\mathbf{a}}\mu_{\mathbf{p}}(\mathbf{q})=\mu_{\mathbf{a}\oplus\mathbf{p}}(\mathbf{q}),$\\
(8) $T_{-\mathbf{a}}\mu_{\mathbf{a}}(\mathbf{q})=1,$\\
(9) $T_{\mathbf{a}}\mu_{\mathbf{p}}(\mathbf{q})=T_{-\mathbf{a}}\mu_{\mathbf{q}}(\mathbf{p})\mu_{\mathbf{q}}(\mathbf{q})\mu_{\mathbf{p}}(\mathbf{0}),$\\
(10) $T_{\mathbf{a}}f(\mathbf{q})=T_{\mathbf{q}}(f_{-}\circ \psi^{\sharp})(\mathbf{a})\mu_{\mathbf{a}}(\mathbf{0})\mu_{\mathbf{q}}(\mathbf{q}),$\\
(11) $T_{\mathbf{b} \oplus \mathbf{a}} f(\mathbf{q})=T_{\mathbf{b}} T_{\mathbf{a}} f^{\sharp}_{\mathbf{b}, \mathbf{a}}(\mathbf{q}),$\\
(12) $T_{-\mathbf{a}}T_{\mathbf{a}} f(\mathbf{q}) = f(\mathbf{q}),$\\
where $f^{\sharp}_{\mathbf{b}, \mathbf{a}}:=f\circ\operatorname{gyr}[\mathbf{b}, \mathbf{a}]$.
\end{prop}
\begin{proof}
Denote $\lambda=\frac{1+\sigma}{2}$.
(1)-(2) follows by the definition of $\mu_{\mathbf{a}}$ in \eqref{def mu}.\\
(3) By definition of $\mu_{\mathbf{a}}$ in \eqref{def mu},
\begin{equation*}\begin{aligned}
\mu_{\mathbf{a}}(\mathbf{q})=\left(\frac{1-|\mathbf{a}|^2}{1-|\mathbf{a}^*\mathbf{q}|^2}\right)^{\lambda}=\left(\frac{1-|\mathbf{q}|^2}{1-|\mathbf{q}^*\mathbf{a}|^2}\right)^{\lambda}\cdot \left(\frac{1-|\mathbf{a}|^2}{1-|\mathbf{q}|^2}\right)^{\lambda}=\mu_{\mathbf{q}}(\mathbf{a})\mu_{\mathbf{a}}(\mathbf{0})\mu_{\mathbf{q}}(\mathbf{q}).
\end{aligned}\end{equation*}
(4) It follows from definitions of $\mu_{\mathbf{a}}$ and $\mathbf{a}\oplus\mathbf{q}$ that
\begin{equation*}\begin{aligned}
\mu_{\mathbf{a}}(\mathbf{a}\oplus\mathbf{q})=&\left(\frac{1-|\mathbf{a}|^2}{|1-\mathbf{a}^*(\mathbf{a}\oplus\mathbf{q})|^2}\right)^{\lambda}=\left(\frac{1-|\mathbf{a}|^2}{|1-\mathbf{a}^*(\mathbf{a}+\mathbf{A}\mathbf{q})(1+\mathbf{a}^*\mathbf{q})^{-1}|^2}\right)^{\lambda}\\
=&\left(\frac{1-|\mathbf{a}|^2}{|1-(|\mathbf{a}|^2+\mathbf{a}^*\mathbf{q})(1+\mathbf{a}^*\mathbf{q})^{-1}|^2}\right)^{\lambda}\\
=&\left(\frac{|1+\mathbf{a}^*\mathbf{q}|^2}{1-|\mathbf{a}|^2}\right)^{\lambda}=\mu_{-\mathbf{a}}(\mathbf{q})^{-1}.\\
\end{aligned}\end{equation*}
(5) By the identity \eqref{module varb} and \eqref{def mu2}, we have
\begin{equation*}\begin{aligned}
\mu_{(-\mathbf{a})\oplus\mathbf{q}}(\mathbf{0})=&\left(1-|(-\mathbf{a})\oplus\mathbf{q}|^2\right)^{\lambda}=\left(1-|\mathbf{q}\oplus(-\mathbf{a})|^2\right)^{\lambda}=\mu_{\mathbf{q}\oplus(-\mathbf{a})}(\mathbf{0}),\\
\mu_{\mathbf{a}}(\mathbf{q})\mu_{\mathbf{q}}(\mathbf{0})=&\left(1-|\varphi_{\mathbf{a}}(\mathbf{q})|^2\right)^{\lambda}=\left(\frac{1-|\varphi_{\mathbf{q}}(\mathbf{a})|^2}{1-|\mathbf{a}|^2}\right)^{\lambda}\cdot(1-|\mathbf{a}|^2)^{\lambda}=\mu_{\mathbf{q}}(\mathbf{a})\mu_{\mathbf{a}}(\mathbf{0}),
\end{aligned}\end{equation*}
and they both equal $\left(1-|\varphi_{\mathbf{a}}(\mathbf{q})|^2\right)^{\lambda}.$\\
(6) We have
\begin{equation*}\begin{aligned}
\mu_{(-\mathbf{a})\oplus\mathbf{q}}((-\mathbf{a})\oplus\mathbf{q})=\mu_{(-\mathbf{a})\oplus\mathbf{q}}(\mathbf{0})^{-1}=(\mu_{\mathbf{a}}(\mathbf{q})\mu_{\mathbf{q}}(\mathbf{0}))^{-1}=\mu_{\mathbf{a}}(\mathbf{q})^{-1}\mu_{\mathbf{q}}(\mathbf{q})
\end{aligned}\end{equation*}
by $(2)$, $(5)$, (3) respectively.\\
(7) Note that
\begin{equation}\begin{aligned}\label{3.6}
T_{\mathbf{a}}\mu_{\mathbf{p}}(\mathbf{q})=&\mu_{\mathbf{a}}(\mathbf{q})\mu_{\mathbf{p}}((-\mathbf{a})\oplus \mathbf{q})=\mu_{\mathbf{a}}(\mathbf{q})\mu_{\mathbf{p}}(\varphi_{-\mathbf{a}}(-\mathbf{q}))\\
=&\left(\frac{1-|\varphi_{\mathbf{a}}(\mathbf{q})|^2}{1-|\mathbf{q}|^2}\right)^{\lambda}\cdot\left(\frac{1-|\varphi_{\mathbf{p}}(\varphi_{-\mathbf{a}}(-\mathbf{q}))|^2}{1-|\varphi_{-\mathbf{a}}(-\mathbf{q})|^2}\right)^{\lambda}\\
=&\left(\frac{1-|\varphi_{\mathbf{p}}(\varphi_{-\mathbf{a}}(-\mathbf{q}))|^2}{1-|\mathbf{q}|^2}\right)^{\lambda}
\end{aligned}\end{equation}
by \eqref{def mu2}. Since
\begin{equation}\begin{aligned}\label{3.8}
\varphi_{\mathbf{p}}(\varphi_{-\mathbf{a}}(-\mathbf{q}))=&\mathbf{p}\oplus(-((-\mathbf{a})\oplus\mathbf{q}))=\mathbf{p}\oplus(\mathbf{a}\oplus(-\mathbf{q}))\\
=&(\mathbf{p}\oplus\mathbf{a})\oplus\mathrm{gyr}[\mathbf{p},\mathbf{a}](-\mathbf{q})=\mathrm{gyr}[\mathbf{p},\mathbf{a}](\mathbf{a}\oplus\mathbf{p})\oplus\mathrm{gyr}[\mathbf{p},\mathbf{a}](-\mathbf{q})\\
=&\mathrm{gyr}[\mathbf{p},\mathbf{a}]\left((\mathbf{a}\oplus\mathbf{p})\oplus(-\mathbf{q})\right)=\mathrm{gyr}[\mathbf{p},\mathbf{a}]\varphi_{\mathbf{a}\oplus\mathbf{p}}(\mathbf{q}),
\end{aligned}\end{equation}
by using gyrocommutativity and gyroautomorphism, we get
\begin{equation}\begin{aligned}\label{mu int}
T_{\mathbf{a}}\mu_{\mathbf{p}}(\mathbf{q})=&\left(\frac{1-|\varphi_{\mathbf{a}\oplus\mathbf{p}}(\mathbf{q})|^2}{1-|\mathbf{q}|^2}\right)^{\lambda}=\mu_{\mathbf{a}\oplus\mathbf{p}}(\mathbf{q})\\
\end{aligned}\end{equation}
by using the isometry of the gyration operator and \eqref{def mu2}.\\
(8) By (7), we get $T_{-\mathbf{a}}\mu_{\mathbf{a}}(\mathbf{q})=\mu_{(-\mathbf{a})\oplus\mathbf{a}}(\mathbf{q})=\mu_{\mathbf{0}}(\mathbf{q})=1.$\\
(9) By applying \eqref{module varb} to the first identity in \eqref{mu int}, we get
\begin{equation*}\begin{aligned}
T_{\mathbf{a}}\mu_{\mathbf{p}}(\mathbf{q})
=&\left(\frac{1-|\mathbf{a}|^2}{|1+\mathbf{a}^*\mathbf{p}|^2}\cdot\frac{1-|\mathbf{q}|^2}{|1-\mathbf{q}^*(\mathbf{a}\oplus\mathbf{p})|^2}\cdot\frac{1}{1-|\mathbf{q}|^2}\cdot(1-|\mathbf{p}|^2)\right)^{\lambda}\\
=&\mu_{-\mathbf{a}}(\mathbf{p})\mu_{\mathbf{q}}(\mathbf{a}\oplus\mathbf{p})\mu_{\mathbf{q}}(\mathbf{q})\mu_{\mathbf{p}}(\mathbf{0})
=T_{-\mathbf{a}}\mu_{\mathbf{q}}(\mathbf{p})\mu_{\mathbf{q}}(\mathbf{q})\mu_{\mathbf{p}}(\mathbf{0}).
\end{aligned}\end{equation*}
(10) We have
\begin{equation*}\begin{aligned}
T_{\mathbf{q}}[f_{-}\circ\psi^{\sharp}_{\mathbf{q}}](\mathbf{a})\mu_{\mathbf{a}}(\mathbf{0})\mu_{\mathbf{q}}(\mathbf{q})=&f\bigg(-\mathrm{gyr}[\mathbf{q},(-\mathbf{q})\oplus\mathbf{a}]((-\mathbf{q})\oplus\mathbf{a})\bigg)\mu_{\mathbf{q}}(\mathbf{a})\mu_{\mathbf{a}}(\mathbf{0})\mu_{\mathbf{q}}(\mathbf{q})\\
=&f\bigg(-\mathrm{gyr}[\mathbf{a},-\mathbf{q}]((-\mathbf{q})\oplus\mathbf{a})\bigg)\mu_{\mathbf{a}}(\mathbf{q})\\
=&f(-(\mathbf{a}\oplus(-\mathbf{q}))\mu_{\mathbf{a}}(\mathbf{q})\\
=&f((-\mathbf{a})\oplus\mathbf{q})\mu_{\mathbf{a}}(\mathbf{q})=T_{\mathbf{a}}f(\mathbf{q}),
\end{aligned}\end{equation*}
by definition \eqref{def T}, (3), where we apply the special loop property in Proposition \ref{aboplus} with $\mathbf{c}=\mathbf{q}$, $\mathbf{d}=(-\mathbf{q})\oplus \mathbf{a}$ in the second identity, and gyrocommutativity in the third identity.\\
(11) Applying definition \eqref{def T} twice, we get
\begin{equation}\begin{aligned}\label{11}
T_{\mathbf{b}} (T_{\mathbf{a}} f)(\mathbf{q})
=& f\bigg((-\mathbf{a}) \oplus \bigg((-\mathbf{b}) \oplus \mathbf{q}\bigg)\bigg)  \mu_{\mathbf{a}}((-\mathbf{b})\oplus\mathbf{q})\mu_{\mathbf{b}}(\mathbf{q})\\
=& f\bigg(\operatorname{gyr}[-\mathbf{a}, -\mathbf{b}]\bigg[\bigg((-\mathbf{b}) \oplus (-\mathbf{a})\bigg) \oplus \mathbf{q}\bigg]\bigg)T_{\mathbf{b}} \mu_{\mathbf{a}}(\mathbf{q})\\
=& f\bigg(\operatorname{gyr}[\mathbf{a}, \mathbf{b}]\bigg[\bigg(-(\mathbf{b} \oplus \mathbf{a})\bigg) \oplus \mathbf{q}\bigg]\bigg) \mu_{\mathbf{b} \oplus \mathbf{a}}(\mathbf{q})\\
=& T_{\mathbf{b} \oplus \mathbf{a}} \bigg( f\circ\operatorname{gyr}[\mathbf{a}, \mathbf{b}]\bigg)(\mathbf{q}),
\end{aligned}\end{equation}
where we used Lemma \ref{cancellation} (i) in the second idnetity and (7) in the third identity. The identity is proved if we substitute $f$ by $f(\operatorname{gyr}[\mathbf{b}, \mathbf{a}]\cdot)$, since $\operatorname{gyr}[\mathbf{b}, \mathbf{a}]\cdot\operatorname{gyr}[\mathbf{a}, \mathbf{b}]={\rm id}$ by Lemma \ref{cancellation} (ii).\\
(12) follows from \eqref{11} by setting \( \mathbf{b} = -\mathbf{a} \), since \( (-\mathbf{a}) \oplus \mathbf{a} = \mathbf{0} \) and \(\operatorname{gyr}[\mathbf{a}, -\mathbf{a}] = id\).
\end{proof}
Recall that the Jacobian $J_R(\varphi_{\mathbf{a}})$ of the smooth M\"{o}bius transformation $\varphi_{\mathbf{a}}$ from $B^{4n}$ to itself is
\begin{equation}\begin{aligned}\label{jacobian}
J_R(\varphi_{\mathbf{a}})(\mathbf{q})=\bigg(\frac{1-|\varphi_{\mathbf{a}}(\mathbf{q})|^2}{1-|\mathbf{q}|^2}\bigg)^{2n+2},
\end{aligned}\end{equation}
for $\mathbf{a}\in B^{4n}$, by \cite[Lemma 5.1]{Wang}.
The adjoint operator of $T_{\mathbf{a}}$ is $T_{-\mathbf{a}}$ by the following proposition.
\begin{prop}\label{prop inner}
For \( f, g \in L^2(B^{4n}, dV_{\sigma}) \) and \( \mathbf{a} \in B^{4n} \), $\sigma\in \mathbb{R}$, we have
\begin{equation}\begin{aligned}\label{eq inner}
\int_{B^{4n}} T_{\mathbf{a}} f(\mathbf{p}) \, g(\mathbf{p}) \, dV_{\sigma}(\mathbf{p})=\int_{B^{4n}} f(\mathbf{p}) \, T_{-\mathbf{a}} g(\mathbf{p}) \, dV_{\sigma}(\mathbf{p}).
\end{aligned}\end{equation}
\end{prop}

\begin{proof}
After changing variables \(\mathbf{p} =\mathbf{a} \oplus \mathbf{w}\), the measure changes as
\begin{equation}\begin{aligned}\label{dVa+w}
dV_{\sigma}(\mathbf{a} \oplus \mathbf{w})=& \left( 1 - |\varphi_{\mathbf{a}}(-\mathbf{w})|^2 \right)^{\sigma-2n-1}d\left(\varphi_{\mathbf{a}}\left(-\mathbf{w}\right)\right)\\
=&\left( 1 - |\varphi_{\mathbf{a}}(-\mathbf{w})|^2 \right)^{\sigma-2n-1}\left( \frac{1 - |\varphi_{\mathbf{a}}(-\mathbf{w})|^2}{1 - |\mathbf{w}|^2} \right)^{2n+2}
   d\mathbf{w} \\
=& \mu_{-\mathbf{a}}(\mathbf{w})^2 \, dV_{\sigma}(\mathbf{w}),
\end{aligned}\end{equation}
by \eqref{jacobian} and definition of $\mu_{-\mathbf{a}}$ in \eqref{def mu2}.
Then, we have
\begin{equation}\begin{aligned}
LHS\; {\rm of} \; \eqref{eq inner}&= \int_{B^{4n}} f((-\mathbf{a}) \oplus \mathbf{p}) \, \mu_{\mathbf{a}}(\mathbf{p}) \, g(\mathbf{p}) \, dV_{\sigma}(\mathbf{p})\\
 &= \int_{B^{4n}} f(\mathbf{w}) \, \mu_{\mathbf{a}}(\mathbf{a} \oplus \mathbf{w}) \,
  g(\mathbf{a} \oplus \mathbf{w})\mu_{-\mathbf{a}}(\mathbf{w})^2 \, dV_{\sigma}(\mathbf{w}) \\
&= \int_{B^{4n}} f(\mathbf{w}) \,
  g(\mathbf{a} \oplus \mathbf{w}) \,
   T_{-\mathbf{a}} \mu_{\mathbf{a}}(\mathbf{w}) \, \mu_{-\mathbf{a}}(\mathbf{w}) \,
   dV_{\sigma}(\mathbf{w}) \\
&= \int_{B^{4n}} f(\mathbf{w}) \, T_{-\mathbf{a}} g(\mathbf{w}) \, dV_{\sigma}(\mathbf{w}),
\end{aligned}\end{equation}
by changing variable \(\mathbf{p} =\mathbf{a} \oplus \mathbf{w}\) and using Proposition \ref{lem mu} (8).
\end{proof}
Substituting $f$ by $f_{-}$ in \eqref{eq inner} and taking coordinates change $\mathbf{p}\rightarrow-\mathbf{p}$, we get
\begin{equation}\begin{aligned}\label{eq inner 1}
\int_{B^{4n}} T_{\mathbf{a}}[f_{-}](\mathbf{p})g(\mathbf{p}) dV_{\sigma}(\mathbf{p})=&\int_{B^{4n}} f_{-}(\mathbf{p})T_{-\mathbf{a}}g(\mathbf{p}) dV_{\sigma}(\mathbf{p})\\
=&\int_{B^{4n}} f(-\mathbf{p})\mu_{-\mathbf{a}}(\mathbf{p})g(\mathbf{a}\oplus\mathbf{p}) dV_{\sigma}(\mathbf{p})\\
=&\int_{B^{4n}} f(\mathbf{p})\mu_{-\mathbf{a}}(-\mathbf{p})g(\mathbf{a}\oplus(-\mathbf{p})) dV_{\sigma}(\mathbf{p})\\
=&\int_{B^{4n}} f(\mathbf{p}) T_{\mathbf{a}}[g_{-}](\mathbf{p}) dV_{\sigma}(\mathbf{p}).
\end{aligned}\end{equation}
\begin{proof}[Proof of Theorem \ref{Young}]
For \( p = 1 \) and \( r = q \), by using Minkowski inequality, we get
\begin{equation}
\begin{aligned}\label{3.12}
\|f * g\|_q&= \left( \int_{B^{4n}} \left| \int_{B^{4n}} f(\mathbf{q}) T_{\mathbf{p}} [g_{-}](\mathbf{q}) \mu_{\mathbf{p}}(\mathbf{p}) \, dV_{\sigma}(\mathbf{q}) \right|^q dV_{\sigma}(\mathbf{p}) \right)^{1/q} \\
&\leq \int_{B^{4n}} |f(\mathbf{q})| \left( \int_{B^{4n}} \bigg|[T_{\mathbf{p}} [g_{-}](\mathbf{q}) \mu_{\mathbf{p}}(\mathbf{p})\bigg|^q dV_{\sigma}(\mathbf{p}) \right)^{1/q} dV_{\sigma}(\mathbf{q})\\
&\leq \int_{B^{4n}} |f(\mathbf{q})| \left( \int_{B^{4n}} \bigg|g(\mathbf{p} \oplus (-\mathbf{q})) \mu_{\mathbf{p}}(\mathbf{q}) \mu_{\mathbf{p}}(\mathbf{p})\bigg|^q \, dV_{\sigma}(\mathbf{p}) \right)^{1/q} dV_{\sigma}(\mathbf{q}).
\end{aligned}
\end{equation}
Note that
\begin{equation}\begin{aligned}\label{reformation lem}
\mu_{\mathbf{p}}(\mathbf{q}) \mu_{\mathbf{p}}(\mathbf{p})=\mu_{\mathbf{q}}(\mathbf{p}) \mu_{\mathbf{q}}(\mathbf{q}) =|1-\mathbf{q}^*\mathbf{p}|^{-\sigma-1}\leq 2^{-\sigma-1},
\end{aligned}\end{equation}
by Proposition \ref{lem mu} (2)-(3), since $-\sigma-1>0$.
Then,
\begin{equation*}
\begin{aligned}
\|f * g\|_q&\leq \int_{B^{4n}} |f(\mathbf{q})| \left( \int_{B^{4n}} \bigg|g\bigg(\operatorname{gyr}[\mathbf{p}, -\mathbf{q}]\left((-\mathbf{q}) \oplus \mathbf{p}\right)\bigg) \mu_{\mathbf{q}}(\mathbf{p}) \mu_{\mathbf{q}}(\mathbf{q})\bigg|^q \, dV_{\sigma}(\mathbf{p}) \right)^{1/q} dV_{\sigma}(\mathbf{q}) \\
&= \int_{B^{4n}} \bigg|f(\mathbf{q})| \left( \int_{B^{4n}} \bigg|g(\operatorname{gyr}[\mathbf{q} \oplus \mathbf{z}, -\mathbf{q}]\mathbf{z}) \mu_{\mathbf{q}}(\mathbf{q} \oplus \mathbf{z}) \mu_{\mathbf{q}}(\mathbf{q})\bigg|^q \mu_{\mathbf{q}}(-\mathbf{z})^2 \, dV_{\sigma}(\mathbf{z}) \right)^{1/q} dV_{\sigma}(\mathbf{q}) \\
&\leq 2^{-\sigma-1} \int_{B^{4n}} |f(\mathbf{q})| \left( \int_{B^{4n}} \bigg|g(\operatorname{gyr}[\mathbf{q}, \mathbf{z}]\mathbf{z})\bigg|^q \, dV_{\sigma}(\mathbf{z}) \right)^{1/q} dV_{\sigma}(\mathbf{q})  \\
&= 2^{-\sigma-1} \|f\|_1 \|\tilde{g}\|_q,
\end{aligned}
\end{equation*}
by changing variables $\mathbf{p} = \mathbf{q} \oplus \mathbf{z},$ and using \eqref{dVa+w}, \eqref{reformation lem}, where in the second inequality we use the special loop property Proposition \ref{aboplus} with $\mathbf{c}=\mathbf{q}$ and $\mathbf{d}=\mathbf{z}$, i.e.
\begin{equation*}
\begin{aligned}
\operatorname{gyr}[\mathbf{q} \oplus \mathbf{z}, -\mathbf{q}]\mathbf{z}=\operatorname{gyr}[\mathbf{q}, \mathbf{z}]\mathbf{z},
\end{aligned}
\end{equation*}
and
\begin{equation}
\begin{aligned}\label{3.15}
|\mu_{\mathbf{q}}(\mathbf{q} \oplus \mathbf{z}) \mu_{\mathbf{q}}(\mathbf{q})|^q (\mu_{\mathbf{q}}(-\mathbf{z}))^2
=(\mu_{\mathbf{q}}(\mathbf{q} \oplus \mathbf{z}) \mu_{\mathbf{q}}(\mathbf{q}))^{q-1} \mu_{\mathbf{q}}(-\mathbf{z}) \mu_{\mathbf{q}}(\mathbf{q})
\leq2^{-q(\sigma+1)}
\end{aligned}
\end{equation}
by \eqref{reformation lem}. Here $\mu_{\mathbf{q}}(\mathbf{q} \oplus \mathbf{z})\mu_{\mathbf{q}}(-\mathbf{z})=(\mu_{-\mathbf{q}}(\mathbf{z}))^{-1}\mu_{-\mathbf{q}}(\mathbf{z})=1$ by Proposition \ref{lem mu} (4).

  For \( r = \infty \) and \( \frac{1}{p} + \frac{1}{q} = 1 \),
\[
\begin{aligned}
\|f * g\|_{\infty}
&\leq \sup_{\mathbf{p} \in B^{4n}} \int_{B^{4n}} \bigg|f(\mathbf{q}) T_{\mathbf{p}} [g_{-}](\mathbf{q}) \mu_{\mathbf{p}}(\mathbf{p})\bigg| \, dV_{\sigma}(\mathbf{q}) \\
&\leq \sup_{\mathbf{p} \in B^{4n}} \|f\|_p \left( \int_{B^{4n}} \bigg|g(\mathbf{p} \oplus (-\mathbf{q})) \mu_{\mathbf{p}}(\mathbf{q}) \mu_{\mathbf{p}}(\mathbf{p})\bigg|^q \, dV_{\sigma}(\mathbf{q}) \right)^{1/q}  \\
&= \|f\|_p \sup_{\mathbf{p} \in B^{4n}} \left( \int_{B^{4n}} \bigg|g(\mathbf{z}) \mu_{\mathbf{p}}\bigg(\mathbf{p} \oplus (-\mathbf{z})\bigg) \mu_{\mathbf{p}}(\mathbf{p})\bigg|^q \mu_{\mathbf{p}}(\mathbf{z})^2 \, dV_{\sigma}(\mathbf{z}) \right)^{1/q} \\
&\leq 2^{-\sigma-1} \|f\|_p \|g\|_q
\end{aligned}
\]
by changing variables $\mathbf{q} = \mathbf{p} \oplus (-\mathbf{z})$ and using \eqref{3.15} for factors $\mu$ again.

  Now for the general case, let \( g \in L^q(B^{4n}, dV_{\sigma}) \) with \( 1 \leq q \leq \infty \) and consider the linear operator \( T \) defined by \( T_g(f) = f * g \). Then, we have proved that $T_g$ is bounded from $L^1$ to $L^{q}$ and $L^p \to L^\infty$ with
\[
\|T_g(f)\|_q \leq 2^{-\sigma-1} \|\tilde{g}\|_q \|f\|_1,\qquad \|T_g(f)\|_{\infty} \leq 2^{-\sigma-1} \|g\|_q \|f\|_p,
\]
where \( 1/p + 1/q = 1 \). Applying the Riesz-Thorin interpolation theorem, we get $T_g$ is bounded from $L^p$ to $L^{r}$, with
\[
\|T_g(f)\|_r \leq 2^{-1-\sigma} \|\tilde{g}\|_q^{1-s} \|g\|_q^s \|f\|_p
\]
with \( \frac{1}{p} + \frac{1}{q} = 1 + \frac{1}{r} \) and \( s = 1 - \frac{q}{r} \). The Theorem is proved.
\end{proof}

Another important property of the $\sigma$-convolution is its $\sigma$-translation invariance. For the real case, it was proved by Ferreira in \cite{Ferreira}.

\begin{thm}
The \( \sigma \)-convolution is $\sigma$-translation invariant in the following sense,
\begin{equation}\begin{aligned}
T_{\mathbf{a}}(f * g)(\mathbf{p}) = (T_{\mathbf{a}} f) * g^{\sharp}_{\mathbf{a}, \mathbf{p}} (\mathbf{p}).
\end{aligned}\end{equation}
\end{thm}
\begin{proof}
Note that
\[
\begin{aligned}
T_{\mathbf{a}}(f * g)(\mathbf{p})
&= (f * g)((-\mathbf{a}) \oplus \mathbf{p}) \mu_{\mathbf{a}}(\mathbf{p}) \\
&= \mu_{(-\mathbf{a}) \oplus \mathbf{p}}((-\mathbf{a}) \oplus \mathbf{p}) \mu_{\mathbf{a}}(\mathbf{p})\int_{B^{4n}} f(\mathbf{q}) T_{(-\mathbf{a}) \oplus \mathbf{p}} [g_{-}](\mathbf{q}) \, dV_{\sigma}(\mathbf{q}),
\end{aligned}
\]
by definition \eqref{def convo} of the $\sigma$-convolution at $(-\mathbf{a}) \oplus \mathbf{p}$.
Since
$T_{(-\mathbf{a}) \oplus \mathbf{p}} [g_{-}](\mathbf{q})=T_{-\mathbf{a}}T_{\mathbf{p}}[g( -\operatorname{gyr}[-\mathbf{a}, \mathbf{p}]\cdot)](\mathbf{q}) $ by applying Proposition \ref{lem mu} (11) to $T_{(-\mathbf{a}) \oplus \mathbf{p}}$ and Lemma \ref{cancellation} (2), we get
\[
\begin{aligned}
T_{\mathbf{a}}(f * g)(\mathbf{p})
&= \mu_{\mathbf{p}}(\mathbf{p})\int_{B^{4n}} f(\mathbf{q}) T_{-\mathbf{a}} T_{\mathbf{p}} [g(-\operatorname{gyr}[-\mathbf{a}, \mathbf{p}] \cdot)](\mathbf{q}) \, dV_{\sigma}(\mathbf{q}) \\
&= \mu_{\mathbf{p}}(\mathbf{p})\int_{B^{4n}} T_{\mathbf{a}} f(\mathbf{q}) T_{\mathbf{p}} [g(-\operatorname{gyr}[-\mathbf{a}, \mathbf{p}] \cdot)](\mathbf{q}) \, dV_{\sigma}(\mathbf{q})  \\
&= (T_{\mathbf{a}} f) * g(\operatorname{gyr}[-\mathbf{a}, \mathbf{p}] \, \cdot)(\mathbf{p}),
\end{aligned}
\]
by using $\mu_{\mathbf{a}\oplus(-\mathbf{p})}(\mathbf{a} \oplus (-\mathbf{p})) \mu_{\mathbf{a}}(\mathbf{p})=\mu_{\mathbf{p}}(\mathbf{p})$, Proposition \ref{prop inner}, \eqref{reformation lem} and definition \eqref{def convo}.
\end{proof}

\begin{prop}\label{convo associate}
If \( f, g, h \in L^1(B^{4n}, dV_{\sigma}) \), then \eqref{as eq} holds.
\end{prop}
\begin{proof}
By definition \eqref{def mu} and \eqref{def convo}, we have
\begin{equation}
\begin{aligned}\label{f convo 1}
(f * (g * h))(\mathbf{a})&= \mu_{\mathbf{a}}(\mathbf{a})\int_{B^{4n}} f(\mathbf{p}) T_{\mathbf{a}}[(g * h)(-\cdot)](\mathbf{p})dV_{\sigma}(\mathbf{p})  \\
&= \mu_{\mathbf{a}}(\mathbf{a})\int_{B^{4n}} f(\mathbf{p}) (g * h)(\mathbf{a} \oplus (-\mathbf{p})) \mu_{\mathbf{a}}(\mathbf{p}) \, dV_{\sigma}(\mathbf{p})\\
&= \mu_{\mathbf{a}}(\mathbf{a})\int_{B^{4n}} f(\mathbf{p})dV_{\sigma}(\mathbf{p}) \int_{B^{4n}} g(\mathbf{q}) T_{\mathbf{a}\oplus(-\mathbf{p})} [h_{-}](\mathbf{q}) \mu_{\mathbf{p}}(\mathbf{p})  dV_{\sigma}(\mathbf{q}),  \\
\end{aligned}
\end{equation}
by definition \eqref{def convo} of $\sigma$-convolution, Proposition \ref{lem mu} (1),(6) and Proposition \ref{prop inner}. Note that by \eqref{eq inner 1},
\[
\begin{aligned}
\int_{B^{4n}}g(\mathbf{q}) T_{\mathbf{a}\oplus(-\mathbf{p})} [h_{-}](\mathbf{q})dV_{\sigma}(\mathbf{q})=&\int_{B^{4n}}T_{\mathbf{a}\oplus(-\mathbf{p})}[g_{-}](\mathbf{q}) h(\mathbf{q})dV_{\sigma}(\mathbf{q})\\
=&\int_{B^{4n}} T_{\mathbf{a}}T_{-\mathbf{p}} [g(-\operatorname{gyr}[\mathbf{a}, -\mathbf{p}]\cdot)](\mathbf{q})h(\mathbf{q})dV_{\sigma}(\mathbf{q})\\
=&\int_{B^{4n}} T_{-\mathbf{p}} [g(-\operatorname{gyr}[\mathbf{a}, -\mathbf{p}]\cdot)](\mathbf{q})T_{-\mathbf{a}}h(\mathbf{q})dV_{\sigma}(\mathbf{q})\\
=&\int_{B^{4n}}  g(-\operatorname{gyr}[\mathbf{a}, -\mathbf{p}](\mathbf{p}\oplus\mathbf{q}))\mu_{-\mathbf{p}}(\mathbf{q})h(\mathbf{a} \oplus \mathbf{q}) \mu_{-\mathbf{a}}(\mathbf{q})dV_{\sigma}(\mathbf{q})\\
=&\int_{B^{4n}} T_{\mathbf{p}} [g(\operatorname{gyr}[\mathbf{a}, -\mathbf{p}]\cdot)](-\mathbf{q})T_{\mathbf{a}}[h_{-}](-\mathbf{q})dV_{\sigma}(\mathbf{q})\\
=&\int_{B^{4n}} T_{\mathbf{p}} [g(\operatorname{gyr}[\mathbf{a}, -\mathbf{p}]\cdot)](\mathbf{q})T_{\mathbf{a}}[h_{-}](\mathbf{q})dV_{\sigma}(\mathbf{q}),
\end{aligned}
\]
by using Proposition \ref{lem mu} (11) and changing variables $\mathbf{q} = -\mathbf{q}$ in the last identity. Substituting this identity into \eqref{f convo 1}, we get
\[
\begin{aligned}
(f * &(g * h))(\mathbf{a})= \mu_{\mathbf{a}}(\mathbf{a})\int_{B^{4n}} T_{\mathbf{a}} [h_{-}](\mathbf{q}) dV_{\sigma}(\mathbf{q})\int_{B^{4n}} f(\mathbf{p}) T_{\mathbf{p}} [g(\operatorname{gyr}[\mathbf{a}, -\mathbf{p}]\cdot)](\mathbf{q}) \mu_{\mathbf{p}}(\mathbf{p})dV_{\sigma}(\mathbf{p})   \\
=& \mu_{\mathbf{a}}(\mathbf{a}) \int_{B^{4n}} T_{\mathbf{a}} [h_{-}](\mathbf{q})dV_{\sigma}(\mathbf{q}) \int_{B^{4n}} f(\mathbf{p}) T_{\mathbf{q}} [g(\operatorname{gyr}[\mathbf{a}, -(\mathbf{q} \oplus \cdot)] \operatorname{gyr}[\mathbf{q},\cdot]\cdot)](\mathbf{p}) \mu_{\mathbf{q}}(\mathbf{q}) \, dV_{\sigma}(\mathbf{p})\\
=& \bigg( f * g(\operatorname{gyr}[\mathbf{a}, -\mathbf{p}] \operatorname{gyr}[\mathbf{q}, \cdot]\cdot) \bigg) *_{\mathbf{q}} h \bigg)(\mathbf{a}),
\end{aligned}
\]
where we used Fubini's theorem, and Proposition \ref{lem mu} (2) and (10), and
\[
\begin{aligned}
T_{\mathbf{p}} [g(\operatorname{gyr}[\mathbf{a}, -\mathbf{p}]\cdot)](\mathbf{q})&=  T_{\mathbf{q}} [g(\operatorname{gyr}[\mathbf{a}, -\cdot] \psi_{\mathbf{q}}^{\sharp}(\cdot)](\mathbf{p})\mu_{\mathbf{p}}(\mathbf{0}) \mu_{\mathbf{q}}(\mathbf{q})\\
\end{aligned}
\]
by substituting $f$ by $g(\operatorname{gyr}[\mathbf{a}, -\mathbf{p}]\cdot)$ and $\mathbf{a}$ by $\mathbf{p}$ in Proposition \ref{lem mu} (10).
\end{proof}
\begin{rem}\label{rem conv}
Gyroassociativity \eqref{as eq} is similar to the real case \cite{Liu1}. If $g$ is a radial function, then $\sigma$-convolution becomes associative, i.e. $(f * (g * h))(\mathbf{a}) = ((f * g) * h)(\mathbf{a})$.
\end{rem}

\section{Weighted quaternionic M\"{o}bius Laplacian operators $\triangle_{\sigma}$}
We will need some tools in quaternionic algebra \cite{Wang} to prove quaternionic $\sigma$-translation invariance of $\triangle_{\sigma}$, including embedding $\tau:\mathbb{H}^n\rightarrow\mathbb{C}^{2n\times 2}$ given by
\begin{equation}\begin{aligned}\label{zaa}
\tau(\mathbf{q})=(z_A^{A'})=\left(\begin{matrix}z_1^{0'}&z_1^{1'}\\z_2^{0'}&z_2^{1'}\\\vdots&\vdots\\z_{n+1}^{0'}&z_{n+1}^{1'}\\z_{n+2}^{0'}&z_{n+2}^{1'}\\\vdots&\vdots\end{matrix}\right):=\left(\begin{matrix}x_0+x_1\mathbf{i}&-x_2-x_3\mathbf{i}\\x_4+x_5\mathbf{i}&-x_6-x_7\mathbf{i}\\\vdots&\vdots\\x_2-x_3\mathbf{i}&\;\;x_0-x_1\mathbf{i}\\x_6-x_7\mathbf{i}&\;\;x_4-x_5\mathbf{i}\\\vdots&\vdots\end{matrix}\right),
\end{aligned}\end{equation}\cite{wang5},
if we write $q_{l}=x_{4l}+\mathbf{i}x_{4l+1}+\mathbf{j}x_{4l+2}+\mathbf{k}x_{4l+3}$, $l=0,1,\dots,n-1$, and complex vector fields
\begin{equation}\begin{aligned}\label{nabla}
\left(\nabla_{A}^{A'}\right)&=\frac{1}{2}\left(\begin{matrix}\nabla_1^{0'}& \nabla_{2}^{0'}&\\\vdots&\vdots\\ \nabla_{n+1}^{1'}&\nabla_{n+2}^{1'}\\\vdots&\vdots&\end{matrix}\right)
&:=\frac{1}{2}\left(\begin{matrix}\partial_{x_0}-\mathbf{i}\partial_{x_1}& -\partial_{x_2}+\mathbf{i}\partial_{x_3}\\\vdots&\vdots\\ \partial_{x_2}+\mathbf{i}\partial_{x_3}&\;\;\partial_{x_0}+\mathbf{i}\partial_{x_1}\\\vdots&\vdots&\end{matrix}\right),
\end{aligned}\end{equation}
where $A=1,2,\dots,2n$, $A'=0',1'$. The quaternionic structure is encoded in these vector fields.
\begin{lem}\label{lem z}\cite[(3.4)-(3.5)]{Wang}
\begin{align}
z_A^{1'}=&\sum_B-J_A^B\overline{z_B^{0'}},\qquad z_A^{0'}=\sum_BJ_A^B\overline{z_B^{1'}},\label{z jab}\\
\sum_{A,B}z_A^{A'}J_A^Bz_B^{A'}=&0,\qquad\sum_{A,B}z_A^{0'}J_A^Bz_B^{1'}=|\mathbf{q}|^2,\qquad  \sum_{A}z_A^{A'}\nabla_A^{A'}|\mathbf{q}|^2=|\mathbf{q}|^2,
\end{align}
where
\begin{equation}\begin{aligned}\label{def J}
J=\left(J_A^B\right):=\left(\begin{matrix}0&I_n\\-I_n&0\end{matrix}\right),
\end{aligned}\end{equation}
$A'=0',1'$, and $A$ and $B$ are row and column indices, respectively.
\end{lem}
\begin{equation}\begin{aligned}\label{tri0}
\triangle_0=&4\sum_{A,B}\nabla_A^{0'}\nabla_{B}^{1'}J_A^B,\qquad R^{A'}=\sum_{A=1}^{2n}z_A^{A'}\nabla_A^{A'}.
\end{aligned}\end{equation}
Namely, the Euclidean norm equals to symplectic product of $\left(z_{A}^{0'}\right)$ and $\left(z_{A}^{1'}\right)$ and $\triangle_0$ is symplect product of $\left(\nabla_{A}^{0'}\right)$ and $\left(\nabla_{A}^{1'}\right)$.

  To prove the invariance of $\triangle_{\sigma}$, we need to know the derivatives of multiplicative factor $\mu_{\mathbf{a}}$ and quaternionic M\"{o}bius transformations $\varphi_{\mathbf{a}}$.
\begin{lem}\label{lem Dmu}
(1) $\left.\nabla_{A}^{A'}\mu_{\mathbf{a}}(\mathbf{q})\right|_{\mathbf{q}=\mathbf{0}}=-(1+\sigma)\mu_{\mathbf{a}}(\mathbf{0})\overline{a_A^{A'}}$,\\
(2) $\left.\triangle_0 \mu_{\mathbf{a}}(\mathbf{q})\right|_{\mathbf{q}=\mathbf{0}}=(\sigma^2-1)\mu_{\mathbf{a}}(\mathbf{0})|\mathbf{a}|^2$.
\end{lem}

\begin{lem}\label{nabla w0}
Denote $(2n\times 2)$-matrix valued function $(w_A^{A'})(\mathbf{q})=\tau(\varphi_{\mathbf{a}}(\mathbf{q})).$ Then, we have
\begin{equation}\begin{aligned}\label{n1}
\left.\nabla_{A}^{0'}w_C^{0'}(\mathbf{q})\right|_{\mathbf{q}=\mathbf{0}}=\left.\nabla_{A}^{1'}w_C^{1'}(\mathbf{q})\right|_{\mathbf{q}=\mathbf{0}}=-S_{\mathbf{a}}\delta_{CA}+\frac{S_{\mathbf{a}}}{1+S_{\mathbf{a}}}(a_C^{0'}\overline{a_A^{0'}}+a_C^{1'}\overline{a_A^{1'}}),\\
\end{aligned}\end{equation}
\begin{equation}\begin{aligned}\label{vanish d}
\left.\nabla_{A}^{0'}w_B^{1'}(\mathbf{q})\right|_{\mathbf{q}=\mathbf{0}}=\left.\nabla_{A}^{1'}w_B^{0'}(\mathbf{q})\right|_{\mathbf{q}=\mathbf{0}}=0.\quad\;\;\;\qquad\qquad\qquad\qquad\qquad\qquad
\end{aligned}\end{equation}
\end{lem}
The proofs of these lemmas will be given in the Appendix. We also need the following chain rule in terms of variables $z_A^{A'}$ and derivatives $\nabla_A^{A'}$.
\begin{lem}\label{chain rule}\cite[Lemma 3.4]{Wang}
For a smooth function $f$ on $\mathbb{R}^{4n}$ and a smooth transformation $\psi$, we have
\begin{equation}\begin{aligned}\label{eq chain}
\nabla_{A}^{A'}(f\circ \psi)(\mathbf{q})&=\sum_{B,B'}\left(\nabla_{B}^{B'}f\right)\left(\psi(\mathbf{q})\right)\cdot\nabla_{A}^{A'}v_{B}^{B'}(\mathbf{q}),\\
\end{aligned}\end{equation}
where $\tau(\psi(\mathbf{q}))=(v_{B}^{B'}(\mathbf{q}))$ is $(2n\times 2)$-matrix valued function.
\end{lem}
\begin{proof}[Proof of Theorem \ref{thm invariance}]
We claim that
\begin{equation}\begin{aligned}\label{delta sigma}
(\triangle_{\sigma}f)(\mathbf{a})=(\mu_{\mathbf{a}}(\mathbf{0}))^{-1}\triangle(T_{-\mathbf{a}} f)(\mathbf{0})-(\sigma^2-1)f(\mathbf{a}).
\end{aligned}\end{equation}
By using the expression \eqref{tri0} of $\triangle_0$ and the chain rule \eqref{eq chain} twice, we get
\begin{equation*}\begin{aligned}
\triangle_0(T_{-\mathbf{a}} f)(\mathbf{0})=&4\sum_{A,B}J_A^B\nabla_{A}^{0'}\nabla_{B}^{1'}(\mu_{-\mathbf{a}}\cdot (f\circ \varphi_{\mathbf{a}}))(\mathbf{0})\\
=&f(\mathbf{a})(\triangle_0 \mu_{-\mathbf{a}})(\mathbf{0})+\mu_{-\mathbf{a}}(\mathbf{0})\triangle_0 (f\circ \varphi_{\mathbf{a}})(\mathbf{0})+\mathscr{R},\\
=&\mu_{\mathbf{a}}(\mathbf{0})(\sigma^2-1)|\mathbf{a}|^2f(\mathbf{a})+\mu_{\mathbf{a}}(\mathbf{0})\triangle_0 (f\circ \varphi_{\mathbf{a}})(\mathbf{0})\\
&-\mu_{\mathbf{a}}(\mathbf{0})(1+\sigma)(1-|\mathbf{a}|^2)(R^{0'}+R^{1'})f(\mathbf{a})\\
=&\mu_{\mathbf{a}}(\mathbf{0})\left((\triangle f)(\mathbf{a})-(1+\sigma)(1-|\mathbf{a}|^2)(R^{0'}+R^{1'})f(\mathbf{a})+|\mathbf{a}|^2(\sigma^2-1)f(\mathbf{a})\right)\\
=&\mu_{\mathbf{a}}(\mathbf{0}) \bigg( (\triangle_{\sigma} f)(\mathbf{a})-(\sigma^2-1) f(\mathbf{a}) \bigg).
\end{aligned}\end{equation*}
by Lemma \ref{lem Dmu}, where the vanishing of \eqref{vanish d} in Lemma \ref{nabla w0} allows us to write
\begin{equation*}\begin{aligned}
\mathscr{R}:=&4\sum_{A,B,C}J_A^B\bigg[(\nabla_{C}^{0'}f)\circ\varphi_{\mathbf{a}}\cdot\nabla_{A}^{0'} w_{C}^{0'}\cdot\nabla_{B}^{1'}\mu_{-\mathbf{a}}+(\nabla_{C}^{1'} f)\circ\varphi_{\mathbf{a}}\cdot\nabla_{B}^{1'}w_{C}^{1'}\cdot\nabla_{A}^{0'}\mu_{-\mathbf{a}}\bigg](\mathbf{0})\\
=&4\sum_{A,B,C}\bigg[\bigg((\nabla_{C}^{0'}f)(\mathbf{a})\bigg(-(1+\sigma)\mu_{\mathbf{a}}(\mathbf{0})\overline{a_B^{1'}}\bigg)\bigg(-S_{\mathbf{a}}\delta_{CA}+\frac{S_{\mathbf{a}}}{1+S_{\mathbf{a}}}(a_C^{0'}\overline{a_A^{0'}}+a_C^{1'}\overline{a_A^{1'}})\bigg)\bigg.\\
&+\bigg.(\nabla_{C}^{1'}f)(\mathbf{a})\bigg(-(1+\sigma)\mu_{\mathbf{a}}(\mathbf{0})\overline{a_A^{0'}}\bigg)\bigg(-S_{\mathbf{a}}\delta_{CB}+\frac{S_{\mathbf{a}}}{1+S_{\mathbf{a}}}(a_C^{0'}\overline{a_B^{0'}}+a_C^{1'}\overline{a_B^{1'}})\bigg)\bigg]J_A^B\\
=&(1+\sigma)\mu_{\mathbf{a}}(\mathbf{0})(|\mathbf{a}|^2-1)\sum_{C}\bigg(a_C^{0'}\nabla_{C}^{0'}f(\mathbf{a}) +a_{C}^{1'}\nabla_{C}^{1'}f(\mathbf{a})\bigg)\\
=& -\mu_{\mathbf{a}}(\mathbf{0})(1+\sigma)(1-|\mathbf{a}|^2)(R^{0'}+R^{1'})f(\mathbf{a}),
\end{aligned}\end{equation*}
by \eqref{n1}, the expression of $R^{A'}$ in \eqref{tri0} and third and fourth identity in \eqref{z jab}.
The claim is proved.
Now note that
\[
\begin{aligned}
T_{-\mathbf{a}} (T_{\mathbf{b}} f)(\mathbf{q})=(T_{\mathbf{b}}f)( \mathbf{a} \oplus \mathbf{q}) \mu_{-\mathbf{a}}(\mathbf{q})=f((-\mathbf{b}) \oplus (\mathbf{a} \oplus \mathbf{q})) T_{-\mathbf{a}} \mu_{\mathbf{b}}(\mathbf{q}),
\end{aligned}
\]
by using definition \eqref{def T} twice and $\mu_{-\mathbf{a}}(\mathbf{q})\mu_{\mathbf{b}}(\mathbf{a}\oplus \mathbf{q})=T_{-\mathbf{a}} \mu_{\mathbf{b}}(\mathbf{q})$ by definition \eqref{def T} of $\sigma$-translation. Applying the claim \eqref{delta sigma} to $T_{\mathbf{b}} f$, we get
\begin{equation}
\begin{aligned}\label{eq 41}
\triangle_{\sigma}(T_{\mathbf{b}} &f)(\mathbf{a})
= \mu_{\mathbf{a}}(\mathbf{0})^{-1} \Delta(T_{-\mathbf{a}} (T_{\mathbf{b}} f))(\mathbf{0})
   -(\sigma^2-1) T_{\mathbf{b}} f(\mathbf{a}) \\
&= \mu_{\mathbf{a}}(\mathbf{0})^{-1}
   \left. \Delta\Big( f((-\mathbf{b}) \oplus (\mathbf{a} \oplus \mathbf{q})) T_{-\mathbf{a}} \mu_{\mathbf{b}}(\mathbf{q}) \Big) \right|_{\mathbf{q}=\mathbf{0}}
   - (\sigma^2-1) f((-\mathbf{b}) \oplus \mathbf{a}) \mu_{\mathbf{b}}(\mathbf{a}).
\end{aligned}
\end{equation}
Note that
\begin{equation}
\begin{aligned}\label{eq 42}
T_{-\mathbf{a}} \mu_{\mathbf{b}}(\mathbf{q})=& \mu_{(-\mathbf{a}) \oplus \mathbf{b}}(\mathbf{q})=\bigg(\frac{1-|\varphi_{(-\mathbf{a}) \oplus \mathbf{b}}(\mathbf{q})|^2}{1-|\mathbf{q}|^2}\bigg)^{\lambda}=\bigg(\frac{1-|((-\mathbf{a}) \oplus \mathbf{b})\oplus(-\mathbf{q})|^2}{1-|\mathbf{q}|^2}\bigg)^{\lambda}\\
=&\bigg(\frac{1-\bigg|\operatorname{gyr}[\mathbf{b}, -\mathbf{a}]((-\mathbf{a}) \oplus \mathbf{b})\oplus\operatorname{gyr}[\mathbf{b}, -\mathbf{a}](-\mathbf{q})\bigg|^2}{1-|\operatorname{gyr}[\mathbf{b}, -\mathbf{a}]\mathbf{q}|^2}\bigg)^{\lambda}\\
=&\bigg(\frac{1-\bigg|( \mathbf{b}\oplus(-\mathbf{a}))\oplus\operatorname{gyr}[-\mathbf{b},\mathbf{a}](\mathbf{q})\bigg|^2}{1-|\mathbf{q}|^2}\bigg)^{\lambda}\\
=& \mu_{\mathbf{b} \oplus (-\mathbf{a})}(\operatorname{gyr}[-\mathbf{b}, \mathbf{a}]\mathbf{q})
\end{aligned}
\end{equation}
by Proposition \ref{lem mu} (7), definition \eqref{def mu2} of factor $\mu$, definition \eqref{def gyr} of quaternionic gyration operator. Here we use unitarity of the gyration operator and gyroautomorphism in the third identity and use gyrocommutavity in the fourth identity.
Then, substituting \eqref{eq 42} into \eqref{eq 41} and substituting $\mathbf{p}=(-\mathbf{b}) \oplus \mathbf{a}$, we get
\[
\begin{aligned}
\triangle_{\sigma}(T_{\mathbf{b}} f)(\mathbf{a})
&= \mu_{\mathbf{a}}(\mathbf{0})^{-1}
   \left. \Delta\Big( f(\mathbf{p} \oplus \operatorname{gyr}[-\mathbf{b}, \mathbf{a}]\mathbf{q})
   \mu_{-\mathbf{p}}(\operatorname{gyr}[-\mathbf{b}, \mathbf{a}]\mathbf{q}) \Big) \right|_{\mathbf{q}=\mathbf{0}} -(\sigma^2-1)f(\mathbf{p} ) \mu_{\mathbf{b}}(\mathbf{a}) \\
&= \left.\mu_{\mathbf{a}}(\mathbf{0})^{-1}
   \Delta\bigg(f (\mathbf{p}\oplus\mathbf{q})\mu_{-\mathbf{p}}(\mathbf{q})\bigg)\right|_{\mathbf{q}=\mathbf{0}}
  -(\sigma^2-1) f(\mathbf{p}) \mu_{\mathbf{b}}(\mathbf{a}) \\
&= \mu_{\mathbf{a}}(\mathbf{0})^{-1}
   \Delta(T_{-\mathbf{p} } f)(\mathbf{0})
  -(\sigma^2-1) f(\mathbf{p}) \mu_{\mathbf{b}}(\mathbf{a}) \\
&= \mu_{\mathbf{b}}(\mathbf{a})\mu_{\mathbf{p} }(\mathbf{0})^{-1}
   \Delta(T_{-\mathbf{p} } f)(\mathbf{0})-(\sigma^2-1) f(\mathbf{p} ) \mu_{\mathbf{b}}(\mathbf{a}) \\
&= (\triangle_{\sigma} f)((-\mathbf{b}) \oplus \mathbf{a}) \mu_{\mathbf{b}}(\mathbf{a})= T_{\mathbf{b}}(\triangle_{\sigma} f)(\mathbf{a}).
\end{aligned}
\]
by  Proposition \ref{lem mu} (1), (5), i.e.
\begin{equation*}
\begin{aligned}
\mu_{\mathbf{p}}(\mathbf{0})=\mu_{\mathbf{a}\oplus(-\mathbf{b})}(\mathbf{0})=\mu_{(-\mathbf{a})\oplus\mathbf{b}}(\mathbf{0})=\mu_{\mathbf{b}}(\mathbf{a})\mu_{\mathbf{a} }(\mathbf{0}),
\end{aligned}
\end{equation*}
where $|\mathbf{p}|=|\mathbf{a}\oplus(-\mathbf{b})|$ by \eqref{module v1}.
We have used the Sp$(n)$Sp$(1)$-invariance \eqref{u in} of $\Delta$ in the second identity and definition \eqref{ex form} of $\triangle_{\sigma}$ in the fifth identity.
Theorem \ref{thm invariance} is proved.
\end{proof}

\section{left gyro-subgroup of quaternionic M\"{o}bius left gyrogroup}
$K$ is called a {\it left gyro-subgroup} of a left gyrogroup $(G, \oplus)$ if it is a {\rm left gyrogroup} for the operation induced from $G$ and $\operatorname{gyr}[a,b] \in \operatorname{Aut}(K)$ for all $a,b \in K$. For a fixed $\boldsymbol{\zeta} \in S^{4n-1}$, we consider the slice $L_{{\boldsymbol{\zeta}}} := \{\boldsymbol{\zeta}q_0  ;  |q_0| < 1,q_0\in \mathbb{H}\}$ and $(L_{{\boldsymbol{\zeta}}})^\perp:=\{\mathbf{q} \in B^{4n};  \mathbf{q}^*\mathbf{p}  = 0 \;{\rm for} \;\mathbf{p}\in L_{{\boldsymbol{\zeta}}}\}$
\begin{prop}\label{D,L}
$(L_{{\boldsymbol{\zeta}}})^\perp$ and $L_{{\boldsymbol{\zeta}}}$ are left gyro-subgroups of $(B^{4n}, \oplus)$.
\end{prop}
\begin{proof}
It is easy to check that $L_{{\boldsymbol{\zeta}}}$ is a left gyrogroup.
For $\mathbf{a}$, $\mathbf{b}\in (L_{{\boldsymbol{\zeta}}})^\perp$, we must have $\mathbf{a}^*\boldsymbol{\zeta}=\mathbf{b}^*\boldsymbol{\zeta}=0$ by definition \eqref{def a oplus b}, and so
\begin{equation}
\begin{aligned}
(\mathbf{a}\oplus\mathbf{b})^*\boldsymbol{\zeta}=&(1+\mathbf{b}^*\mathbf{a})^{-1}(\mathbf{a}^*+\mathbf{b}^*\mathbf{A})\boldsymbol{\zeta}=0,
\end{aligned}
\end{equation}
since
\begin{equation*}
\begin{aligned}
\mathbf{b}^*\mathbf{A}\boldsymbol{\zeta}=\mathbf{b}^*(1-S_{\mathbf{a}})|\mathbf{a}|^{-2}\mathbf{a}\mathbf{a}^*\boldsymbol{\zeta}+S_{\mathbf{a}}\mathbf{b}^*\boldsymbol{\zeta}=0.
\end{aligned}
\end{equation*}
Similarly, $(\mathbf{b}\oplus\mathbf{a})^*\boldsymbol{\zeta}=0$. Thus, $\mathbf{b}\oplus\mathbf{a},\mathbf{a}\oplus\mathbf{b}\in (L_{{\boldsymbol{\zeta}}})^\perp$. Obviously, we also have $\mathbf{0},-\mathbf{a},-\mathbf{b}\in (L_{{\boldsymbol{\zeta}}})^\perp$. Moreover, for $\mathbf{c}\in (L_{{\boldsymbol{\zeta}}})^\perp$,
\begin{equation}
\begin{aligned}
(\mathrm{gyr}[\mathbf{a}, \mathbf{b}]\mathbf{c})^*\boldsymbol{\zeta}=&(1+\mathbf{b}^*\mathbf{a})^{-1}\mathbf{c}^*\mathbf{M}^*A_{\mathbf{d}}^{-1}\boldsymbol{\zeta}
=(1+\mathbf{b}^*\mathbf{a})^{-1}\frac{S_{\mathbf{a}}S_{\mathbf{b}}}{S_{\mathbf{d}}}\mathbf{c}^*\boldsymbol{\zeta}=0,
\end{aligned}
\end{equation}
where $\mathbf{A}_{\mathbf{d}}^{-1}\boldsymbol{\zeta}=S_{\mathbf{d}}^{-1}\boldsymbol{\zeta}$, since $\mathbf{d}^*\boldsymbol{\zeta}=0$. Similarly, $\mathbf{A}\boldsymbol{\zeta}=S_{\mathbf{a}}\boldsymbol{\zeta}$, $\mathbf{B}\boldsymbol{\zeta}=S_{\mathbf{b}}\boldsymbol{\zeta}$. Therefore, $\mathbf{M}^*\boldsymbol{\zeta}=(\mathbf{b}^*\mathbf{a}+\mathbf{B}\mathbf{A})\boldsymbol{\zeta}=S_{\mathbf{a}}S_{\mathbf{b}}\boldsymbol{\zeta}$,
by Lemma \ref{Matrix} (1), (3). Thus, $\mathrm{gyr}[\mathbf{a}$, $\mathbf{b}]\mathbf{c}\in (L_{{\boldsymbol{\zeta}}})^\perp$. The gyroassociativity and gyrocommutativity of $(L_{{\boldsymbol{\zeta}}}^{\perp}, \oplus)$ follows from those for $(B^{4n},\oplus)$. $(L_{{\boldsymbol{\zeta}}}^{\perp}, \oplus)$ is a left gyro-subgroup of $(B^{4n}, \oplus)$.
\end{proof}
\begin{rem}
For $n=1$, our quaternionic gyration operator is the gyration operator in $B^4$. The left gyro-subgroup $L_{{\boldsymbol{\zeta}}}$ is noncommutative.
\end{rem}
\begin{prop}
If \( (P, \oplus) \) is a left gyro-subgroup of \( (B^{4n}, \oplus) \), then \( (P^\perp, \oplus) \) is also a left gyro-subgroup of \( (B^{4n}, \oplus) \).
\end{prop}
\begin{proof}
By Proposition \ref{D,L}, For $\mathbf{a},\mathbf{b},\mathbf{c}\in (P^\perp, \oplus)$ and $\mathbf{d}\in(P, \oplus)$, we have $\mathbf{0},-\mathbf{a},-\mathbf{b},\mathbf{a}\oplus\mathbf{b},\mathrm{gyr}[\mathbf{a}, \mathbf{b}]\mathbf{c}\in\bigcap_{\mathbf{d}\in P} (L_{\mathbf{d}/|\mathbf{d}|})^\perp\subset (P^\perp, \oplus)$. Consequently, \( (P^\perp, \oplus) \) is a left gyrogroup.
\end{proof}
As in the real case, we also have the M\"{o}bius decompositions of \(B^{4n}\).
\begin{thm}\label{thm decomp}
Let \( P \) be a M\"{o}bius left gyro-subgroup of \( B^{4n} \) and \( P^\perp \) its orthogonal complement in \( B^{4n} \). Then for each \( \mathbf{c} \in B^{4n} \), there exist unique \( \mathbf{a}, \mathbf{v} \in P \) and \( \mathbf{b}, \mathbf{u}\in P^\perp \) such that \( \mathbf{c} = \mathbf{a} \oplus \mathbf{b} \) and \( \mathbf{c} = \mathbf{u} \oplus \mathbf{v} \).
\end{thm}

\begin{proof}
 Given \( \mathbf{c} \in B^{4n} \), there exist \( \mathbf{c}_1 \in P \) and \( \mathbf{c}_2 \in P^\perp \) such that $\mathbf{c}_1^*\mathbf{c}_2=\mathbf{c}_2^*\mathbf{c}_1=0$, $\mathbf{c}=\mathbf{c}_1+\mathbf{c}_2$. If \( \mathbf{c}_1 = 0 \) then we can take \( \mathbf{b} = 0 \) and \( \mathbf{a} = \mathbf{c}_2 \). If \( \mathbf{c}_1 \neq 0 \), take \( \mathbf{a} = \mathbf{c}_1q_1 \in P \) and \( \mathbf{b} = \mathbf{c}_2q_2 \in P^\perp \). Then, the equation
\begin{equation*}
\begin{aligned}
( \mathbf{c}_1q_1)\oplus( \mathbf{c}_2q_2)
=( \mathbf{c}_1q_1+A_{ \mathbf{c}_1q_1} \mathbf{c}_2q_2)(1+( \mathbf{c}_1q_1)^*\mathbf{c}_2q_2)^{-1}
=\mathbf{c}_1q_1+S_{\mathbf{c}_1q_1}\mathbf{c}_2q_2=\mathbf{c}_1 + \mathbf{c}_2=\mathbf{c},
\end{aligned}
\end{equation*}
has a unique solution given by $q_1 = 1$, $q_2 =S_{\mathbf{c}_1q_1}^{-1}$. Here $S_{\mathbf{c}_1q_1}$ is a real number. So the first decomposition exists.

To prove the uniqueness of the decomposition, we suppose that there exist \(\mathbf{a}, \mathbf{e} \in P\) and \(\mathbf{b}, \mathbf{f} \in P^\perp\) such that \(\mathbf{c} = \mathbf{a} \oplus \mathbf{b} = \mathbf{e} \oplus \mathbf{f}\). Then,
\begin{equation*}
\begin{aligned}
\mathbf{b} = (-\mathbf{a}) \oplus (\mathbf{e} \oplus \mathbf{f})=((-\mathbf{a}) \oplus \mathbf{e}) \oplus\mathrm{gyr}[\mathbf{a}, \mathbf{e}]\mathbf{f},
 \end{aligned}
\end{equation*}
 by Lemma \ref{lem cancellation} (1) and gyroassociativity. Since $\mathbf{f} \in P^\perp$ and \(\mathbf{a}, \mathbf{e} \in P\), we have
\begin{equation}
\begin{aligned}\label{gyr perp}
\mathrm{gyr}[\mathbf{a}, \mathbf{e}]\mathbf{f}=&\mathbf{A}_{\mathbf{a}\oplus \mathbf{e}}^{-1}(\mathbf{a}\mathbf{e}^*+\mathbf{A}_{\mathbf{a}}\mathbf{A}_{\mathbf{e}})\mathbf{f}\mathbf{q}_{\mathbf{a},\mathbf{e}}^{-1}=S_{\mathbf{a}\oplus \mathbf{e}}^{-1}S_{\mathbf{a}}S_{\mathbf{e}}\mathbf{f}\mathbf{q}_{\mathbf{a},\mathbf{e}}^{-1}\\
=&\bigg(\frac{(1-|\mathbf{a}|^2)(1-|\mathbf{e}|^2)}{1-|\mathbf{a}\oplus \mathbf{e}|^2}\bigg)^{\frac12}\mathbf{f}q_{\mathbf{a},\mathbf{e}}^{-1}=|q_{\mathbf{a},\mathbf{e}}|\mathbf{f}\mathbf{q}_{\mathbf{a},\mathbf{e}}^{-1}\in L_{\mathbf{f}/|\mathbf{f}|}\subset P^{\perp},
\end{aligned}
\end{equation}
where $\mathbf{A}_{\mathbf{e}}\mathbf{f}=S_{\mathbf{e}}\mathbf{f}$, $\mathbf{A}\mathbf{f}=S_{\mathbf{a}}\mathbf{f}$, $\mathbf{A}_{\mathbf{a}\oplus \mathbf{e}}^{-1}\mathbf{f}=S_{\mathbf{a}\oplus \mathbf{e}}^{-1}\mathbf{f}$ by by Lemma \ref{Matrix} (1), (3). By Lemma \ref{lem cancellation} (2), $(-\mathbf{a}) \oplus \mathbf{e}=\mathbf{b}\ominus\mathrm{gyr}[\mathbf{a}, \mathbf{e}]\mathbf{f}\in P^\perp$.
 Since $P$ is a left gyro-subgroup, we have $(-\mathbf{a}) \oplus \mathbf{e}\in P$. Thus,
$(-\mathbf{a}) \oplus \mathbf{e}\in P\bigcap P^\perp=\{\mathbf{0}\}$, i.e. $\mathbf{a}=\mathbf{e}$. Hence $\mathbf{b}=\mathbf{f}$. Another uniqueness is similar.
\end{proof}

\begin{cor}
Let \( P \) be an arbitrary M\"{o}bius left gyro-subgroup of \( B^{4n} \). Then, the family \( \{\boldsymbol{a} \oplus P : \boldsymbol{a} \in P^\perp\} \) is a disjoint partition of \( B^{4n} \), i.e. $B^{4n} =\bigcup_{\boldsymbol{a} \in P^\perp} (\boldsymbol{a} \oplus P).$
\end{cor}

\begin{proof}
We first prove that this family is disjoint. Let \( \boldsymbol{a}, \boldsymbol{c} \in P^\perp \) with \( \boldsymbol{a} \neq \boldsymbol{c} \) and assume that \( (\boldsymbol{a} \oplus P) \cap (\boldsymbol{c} \oplus P) \neq \emptyset \). Then there exists \( \boldsymbol{f} \in B^{4n} \) such that \( \boldsymbol{f} = \boldsymbol{a} \oplus \boldsymbol{b} \) and \( \boldsymbol{f} = \boldsymbol{c} \oplus \boldsymbol{d} \) for some \( \boldsymbol{b}, \boldsymbol{d} \in P \). By uniqueness in Theorem \ref{thm decomp}, \( \boldsymbol{a} = \boldsymbol{c} \). But this contradicts our assumption. Thus, \( (\boldsymbol{a} \oplus P) \cap (\boldsymbol{c} \oplus P) = \emptyset \) provided \( \boldsymbol{a}, \boldsymbol{c} \in P^\perp \) and \( \boldsymbol{a} \neq \boldsymbol{c} \). Finally, by Theorem \ref{thm decomp}, we have \( B^{4n} = \bigcup_{\boldsymbol{a} \in P^\perp} (\boldsymbol{a} \oplus P) \).
\end{proof}
We can construct \emph{quotient M\"{o}bius left gyrogroups}.
\begin{prop}
The quotient space \((B^{4n} / P, \sim_l)\), endowed with the binary operation \(\oplus\) defined by $(\boldsymbol{a} \oplus P) \oplus (\boldsymbol{d} \oplus P) := (\boldsymbol{a} \oplus \boldsymbol{d}) \oplus P,$ with $\boldsymbol{a}, \boldsymbol{d} \in P^{\perp},$ is a quotient M\"{o}bius left gyrogroup.
\end{prop}
It is easy to check conditions (1)-(5) of left gyrogroup for following quotient space. Here we omit the details.
\appendix
\section{Proof of Lemmas}
The quaternionic structure $\tau$ can be extended to map a quaternionic $(p\times m)-$matrix $\mathcal{M}=\mathbf{a}+\mathbf{j}\mathbf{b}$ to complex $2p\times 2m$-matrix $\tau(\mathcal{M})=\left(\begin{matrix}\mathbf{a}&-\overline{\mathbf{b}}\\\mathbf{b}&\overline{\mathbf{a}}\end{matrix}\right)$.

\begin{lem}\cite[Proposition 2.3]{wang2}\label{tau MN}
(1) $\tau(\mathcal{M}\mathcal{N})=\tau(\mathcal{M})\tau(\mathcal{N})$ for a quaternionic $(p\times m)-$matrix $\mathcal{M}$ and a quaternionic $(m\times l)$-matrix $\mathcal{N}$. In particular, for $\mathbf{q}'=\mathcal{M}\mathbf{q}$ with $\mathbf{q},\mathbf{q}'\in \mathbb{H}^n$ and a quaternionic $(n\times n)$-matrix $\mathcal{M}$, we have $\tau(\mathbf{q}')=\tau(\mathcal{M})\tau(\mathbf{q})$ as complex $(2n\times 2)-$matrices.\\
(2)\; Let $I_n$ be the identity $(n\times n)$-matrix. $\mathcal{M}\in {\rm Sp}(n)$ if and only if $\tau(\mathcal{M})J\tau(\mathcal{M})^t=J$.
(3)\;$\tau(\mathcal{M}^*)=\overline{\tau(\mathcal{M})}^t$.
\end{lem}
\begin{lem}\label{lem nabla}\cite[Lemma 3.1]{wang3}
For $z_{B}^{B'}$ and $\nabla_A^{A'}$ given by \eqref{zaa}-\eqref{nabla}, we have
\begin{equation}\begin{aligned}\label{delta}
\nabla_{A}^{A'}z_{B}^{B'}=\delta_{AB}\delta^{A'B'}.
\end{aligned}\end{equation}
\end{lem}
By definition $z_A^{A'}$ in \eqref{zaa} and $J$ in \eqref{def J}, we can check directly

To prove Lemma \ref{lem Dmu}-\ref{nabla w0} , we need the following 2 lemmas about the derivative of $\varphi_{\mathbf{a}}$.
\begin{lem}\label{lem dvar}\cite[Lemma 3.3, (3.15), Proof of Proposition 1.1]{Wang}
$w(\mathbf{q}):=\tau(\varphi_{\mathbf{a}}(\mathbf{q}))$ is a $(2n\times 2)$-matrix valued function of $\mathbf{q}$ satisfying
\begin{equation}\begin{aligned}\label{rewritten eq 1st de AneqB}
&\nabla_{B}^{B'} w(\mathbf{0})=\left(-S_{\mathbf{a}}+\frac{S_{\mathbf{a}}}{1+S_{\mathbf{a}}}\tau(\mathbf{a}\mathbf{a}^*)\right)E_B^{B'},\\
&\triangle_0 w(\mathbf{0})=4(1-|\mathbf{a}|^2)\tau(\mathbf{a}),
\end{aligned}\end{equation}
where $E_B^{B'}$ be the $(2n\times 2)$-matrix with $1$ only in $(B,B')$-entry and $0$ all other entries.
\end{lem}
\begin{lem}\label{lem w}
Assume as in Lemma \ref{lem dvar}. Then,
\begin{align}
\left.\nabla_{A}^{A'}\left(\frac{1-|w(\mathbf{q})|^2}{1-|\mathbf{q}|^2}\right)\right|_{\mathbf{q}=\mathbf{0}}&=-(1-|\mathbf{a}|^2)\overline{a_A^{A'}},\label{A2}\\
\frac14\left.\triangle_0|w(\mathbf{q})|^2\right|_{\mathbf{q}=\mathbf{0}}&=-2nS_{\mathbf{a}}^2.\label{A22}
\end{align}
\end{lem}
\begin{proof}
Assume $A'=0'$. Note that $\tau(\mathbf{a}\mathbf{a}^*)$ is a complex $2n\times 2n$ matrix, and so
\begin{equation}\begin{aligned}\label{nabla w}
\left.\nabla_{A}^{0'}w_C^{0'}(\mathbf{q})\right|_{\mathbf{q}=\mathbf{0}}=&\left(\left(-S_{\mathbf{a}}+\frac{S_{\mathbf{a}}}{1+S_{\mathbf{a}}}\tau(\mathbf{a}\mathbf{a}^*)\right)E_A^{0'}\right)_C^{0'}\\
=&-S_{\mathbf{a}}\delta_{CA}+\frac{S_{\mathbf{a}}}{1+S_{\mathbf{a}}}\tau(\mathbf{a}\mathbf{a}^*)_{C A}=-S_{\mathbf{a}}\delta_{CA}+\frac{S_{\mathbf{a}}}{1+S_{\mathbf{a}}}\bigg(a_C^{0'}\overline{a_A^{0'}}+a_C^{1'}\overline{a_A^{1'}}\bigg),
\end{aligned}\end{equation}
by \eqref{z jab}, \eqref{eq q0}, Lemma \ref{rewritten eq 1st de AneqB}. By taking derivatives and noticing vanishing
\begin{equation}\begin{aligned}\label{eq q0}
\left.\nabla_{A}^{A'}(1-|\mathbf{q}|^2)\right|_{\mathbf{q}=\mathbf{0}}=0\qquad{\rm and}\qquad \left.\nabla_{A}^{0'}w_C^{1'}(\mathbf{q})\right|_{\mathbf{q}=\mathbf{0}}=0
\end{aligned}\end{equation}
by \eqref{vanish d}, we get
\begin{equation*}\begin{aligned}
\left.\nabla_{A}^{0'}\left(\frac{1-|w(\mathbf{q})|^2}{1-|\mathbf{q}|^2}\right)\right|_{\mathbf{q}=\mathbf{0}}=&-\left.(1-|\mathbf{q}|^2)^{-1}\nabla_{A}^{0'}|w(\mathbf{q})|^2\right|_{\mathbf{q}=\mathbf{0}}=-\left.\sum_{C,D}w_D^{1'}J_C^D\nabla_{A}^{0'}w_C^{0'}\right|_{\mathbf{q}=\mathbf{0}}\\
=&-\sum_{C,D}\bigg(-S_{\mathbf{a}}\delta_{CA}+\frac{S_{\mathbf{a}}}{1+S_{\mathbf{a}}}\bigg(a_C^{0'}\overline{a_A^{0'}}+a_C^{1'}\overline{a_A^{1'}}\bigg)\bigg)a_D^{1'}J_C^D\\
=&-\sum_{C,D}\left(S_{\mathbf{a}}a_D^{1'}J_A^D+\frac{S_{\mathbf{a}}}{1+S_{\mathbf{a}}}a_C^{0'}\overline{a_A^{0'}}a_D^{1'}J_C^D\right)\\
=&-S_{\mathbf{a}}\overline{a_A^{0'}}+\frac{S_{\mathbf{a}}}{1+S_{\mathbf{a}}}|\mathbf{a}|^2\overline{a_A^{0'}}=-(1-|\mathbf{a}|^2)\overline{a_A^{0'}},
\end{aligned}\end{equation*}
by using \eqref{z jab} repeatedly, and \eqref{A2}. \eqref{A22} for $A'=1'$ is similar.
Then, taking derivatives more and noticing vanishing as in \eqref{vanish d}, we get
\begin{equation*}\begin{aligned}
\frac14\left.\triangle_0|w(\mathbf{q})|^2\right|_{\mathbf{q}=\mathbf{0}}
=&\sum_{A,B,C,D}\nabla_{A}^{0'}\nabla_{B}^{1'}J_A^B\bigg(w_C^{0'}w_D^{1'}J_C^D\bigg)\\
=&\left.\frac14\sum_{C,D}J_C^D\bigg(\triangle_0w_C^{0'}\cdot w_{D}^{1'}+\triangle_0w_D^{1'}\cdot w_{C}^{0'}\bigg)\right|_{\mathbf{q}=\mathbf{0}}+\mathscr{A}\\
=&2(1-|\mathbf{a}|^2)\left.\sum_{C,D}a_C^{0'}a_D^{1'}J_C^D\right|_{\mathbf{q}=\mathbf{0}}+\mathscr{A}\\
=&2(1-|\mathbf{a}|^2)|\mathbf{a}|^2+2nS_{\mathbf{a}}^2-2|\mathbf{a}|^2S_{\mathbf{a}}^2=2nS_{\mathbf{a}}^2,
\end{aligned}\end{equation*}
by using Lemma \ref{rewritten eq 1st de AneqB}, \eqref{z jab},
where
\begin{equation*}\begin{aligned}
\mathscr{A}:=&\left.\sum_{A,B,C,D}J_A^B\nabla_{A}^{0'}w_C^{0'}\cdot \nabla_{B}^{1'}w_{D}^{1'}J_C^D\right|_{\mathbf{q}=\mathbf{0}}
=\sum_{A,B,C,D}\bigg(-S_{\mathbf{a}}\delta_{CA}+\frac{S_{\mathbf{a}}}{1+S_{\mathbf{a}}}\bigg(a_C^{0'}\overline{a_A^{0'}}+a_C^{1'}\overline{a_A^{1'}}\bigg)\bigg)J_A^B\\
&\qquad\qquad\qquad\qquad\qquad\qquad\qquad\qquad\qquad\quad\;\cdot\bigg(-S_{\mathbf{a}}\delta_{DB}+\frac{S_{\mathbf{a}}}{1+S_{\mathbf{a}}}\bigg(a_D^{0'}\overline{a_B^{0'}}+a_D^{1'}\overline{a_B^{1'}}\bigg)\bigg)J_C^D\\
=&\sum_{A,B}J_A^BJ_A^BS_{\mathbf{a}}^2-2\frac{S_{\mathbf{a}}^2}{1+S_{\mathbf{a}}}\sum_{A,C,D}\bigg(a_C^{0'}\overline{a_A^{0'}}+a_C^{1'}\overline{a_A^{1'}}\bigg)J_A^DJ_C^D\\
&+\bigg(\frac{S_{\mathbf{a}}}{1+S_{\mathbf{a}}}\bigg)^2\sum_{A,B,C,D}J_A^BJ_C^D\bigg(a_C^{1'}\overline{a_A^{1'}}a_D^{0'}\overline{a_B^{0'}}+a_D^{1'}\overline{a_B^{1'}}a_C^{0'}\overline{a_A^{0'}}\bigg)\\
=&2nS_{\mathbf{a}}^2-2\frac{S_{\mathbf{a}}^2}{1+S_{\mathbf{a}}}\sum_{C,D}\bigg(a_C^{0'}a_D^{1'}-a_C^{1'}a_D^{0'}\bigg)J_C^D+2\bigg(\frac{S_{\mathbf{a}}}{1+S_{\mathbf{a}}}\bigg)^2|\mathbf{a}|^4\\
=&2nS_{\mathbf{a}}^2-2\frac{S_{\mathbf{a}}^2}{1+S_{\mathbf{a}}}\cdot 2|\mathbf{a}|^2+2\bigg(\frac{S_{\mathbf{a}}}{1+S_{\mathbf{a}}}\bigg)^2|\mathbf{a}|^4=2nS_{\mathbf{a}}^2-2|\mathbf{a}|^2S_{\mathbf{a}}^2,
\end{aligned}\end{equation*}
by the skew symmetry of $J$, $J^2=-I_{2n}$ and Lemma \ref{lem z} repeatedly.
\end{proof}
\begin{proof}[Proof of Lemma \ref{lem Dmu}]
Note that
\begin{equation}\begin{aligned}
\left.\nabla_{A}^{A'}\mu_{\mathbf{a}}(\mathbf{q})\right|_{\mathbf{q}=\mathbf{0}}=&\left.\nabla_{A}^{A'}\left(\frac{1-|w(\mathbf{q})|^2}{1-|\mathbf{q}|^2}\right)^{\frac{1+\sigma}{2}}\right|_{\mathbf{q}=\mathbf{0}}\\
=&\left.\frac{1+\sigma}{2}\left(\frac{1-|w(\mathbf{q})|^2}{1-|\mathbf{q}|^2}\right)^{\frac{\sigma-1}{2}}\nabla_{A}^{A'}\frac{1-|w(\mathbf{q})|^2}{1-|\mathbf{q}|^2}\right|_{\mathbf{q}=\mathbf{0}}=-(1+\sigma)\mu_{\mathbf{a}}(\mathbf{0})\overline{a_A^{A'}}\\
\end{aligned}\end{equation}
by $w(\mathbf{0})=\tau(\mathbf{a})$, Lemma \ref{lem w}, and vanishing \eqref{A2} of derivatives. (1) is proved.
Then,
\begin{equation*}\begin{aligned}
&\left.\frac14\triangle_0 \mu_{\mathbf{a}}(\mathbf{q})\right|_{\mathbf{q}=\mathbf{0}}=\left.\sum_{A,B}J_A^B\nabla_A^{0'}\nabla_B^{1'}\left(\frac{1-|w(\mathbf{q})|^2}{1-|\mathbf{q}|^2}\right)^{\frac{\sigma+1}{2}}\right|_{\mathbf{q}=\mathbf{0}}\\
&\quad\quad\quad=\frac{1+\sigma}{2}\left(\frac{1-|w(\mathbf{q})|^2}{1-|\mathbf{q}|^2}\right)^{\frac{\sigma-1}{2}}\bigg(\frac{\triangle_0(1-|w(\mathbf{q})|^2)}{4(1-|\mathbf{q}|^2)}\bigg.-\left.\left.\frac{1-|w(\mathbf{q})|^2}{4(1-|\mathbf{q}|^2)^2}\triangle_0(1-|\mathbf{q}|^2)\right)\right|_{\mathbf{q}=\mathbf{0}}+\mathscr{B}\\
&\quad\quad\quad=\frac{1+\sigma}{2}\left(1-|\mathbf{a}|^2\right)^{\frac{\sigma-1}{2}}\bigg(2nS_{\mathbf{a}}^2-2n(1-|\mathbf{a}|^2)\bigg)+(1+\sigma)(\sigma-1)\mu_{\mathbf{a}}(\mathbf{0})|\mathbf{a}|^2\\
&\quad\quad\quad=(1+\sigma)(\sigma-1)\mu_{\mathbf{a}}(\mathbf{0})|\mathbf{a}|^2
\end{aligned}\end{equation*}
by Lemma \ref{lem w} and
\begin{equation*}\begin{aligned}
\frac14\triangle_0(1-|\mathbf{q}|^2)=&\sum_{A,B,C,D}\nabla_{A}^{0'}\nabla_{B}^{1'}J_A^B\bigg(z_C^{0'}z_D^{1'}J_C^D\bigg)=\sum_{A,B}J_A^BJ_A^B=2n,
\end{aligned}\end{equation*}
where
\begin{equation*}\begin{aligned}
\mathscr{B}:=&\left.\frac{(1+\sigma)(\sigma-1)}{4}\sum_{A,B}J_A^B\left(\frac{1-|w(\mathbf{q})|^2}{1-|\mathbf{q}|^2}\right)^{\frac{\sigma-3}{2}}\nabla_{A}^{0'}\left(\frac{1-|w(\mathbf{q})|^2}{1-|\mathbf{q}|^2}\right)\nabla_{B}^{1'}\left(\frac{1-|w(\mathbf{q})|^2}{1-|\mathbf{q}|^2}\right)\right|_{\mathbf{q}=\mathbf{0}}\\
=&\frac{(1+\sigma)(\sigma-1)}{4}\mu_{\mathbf{a}}(\mathbf{0})\cdot4\sum_{A,B}\overline{a_A^{0'}}\overline{a_B^{1'}}J_A^B=(1+\sigma)(\sigma-1)\mu_{\mathbf{a}}(\mathbf{0})|\mathbf{a}|^2.
\end{aligned}\end{equation*}
The Lemma is proved.
\end{proof}

\end{document}